\documentclass[leqno,12pt]{article}
\usepackage[hang]{subfigure}
\usepackage{color}
\usepackage{multirow}
\usepackage{latexsym}
\usepackage{float}
\usepackage{xcolor}
\usepackage{hyperref}

\makeatletter
\newsavebox\myboxA
\newsavebox\myboxB
\newlength\mylenA

\newcommand*\xoverline[2][0.75]{%
    \sbox{\myboxA}{$\m@th#2$}%
    \setbox\myboxB\null
    \ht\myboxB=\ht\myboxA%
    \dp\myboxB=\dp\myboxA%
    \wd\myboxB=#1\wd\myboxA
    \sbox\myboxB{$\m@th\overline{\copy\myboxB}$}
    \setlength\mylenA{\the\wd\myboxA}
    \addtolength\mylenA{-\the\wd\myboxB}%
    \ifdim\wd\myboxB<\wd\myboxA%
       \rlap{\hskip 0.5\mylenA\usebox\myboxB}{\usebox\myboxA}%
    \else
        \hskip -0.5\mylenA\rlap{\usebox\myboxA}{\hskip 0.5\mylenA\usebox\myboxB}%
    \fi}
\makeatother

\newcommand{\var}{\mbox{Var}}
\newcommand{\E}{\mathbb{E}} 
\newcommand{\cov}{\mbox{Cov}} 

\newcommand{\R}{\mathbb R}

\newcommand{\N}{\mathbb{N}}

\newcommand{\F}{\cal F}

\def\P{\mathbb{P}}

\usepackage{verbatim}
\usepackage{amsmath,stackrel}
\usepackage{amssymb}
\usepackage{ntheorem}
\usepackage[ansinew]{inputenc}
\usepackage{graphicx}
\usepackage{epsfig}
\usepackage{dsfont}
\usepackage{bm}
\usepackage{bbm}

\usepackage[authoryear]{natbib}
\newtheorem{theorem}{Theorem}[section]

\newtheorem{example}{Example}[section]

\newtheorem{assumption}{Assumption}[section]
\newtheorem{rem}{Remark}[section]

\begin{document}

\title{Efficient tests for bio-equivalence in functional data}

\author{
{\small Holger Dette, Kevin Kokot} \\
{\small Ruhr-Universit\"at Bochum} \\
{\small Fakult\"at f\"ur Mathematik}\\
{\small Bochum, Germany} \\
{\small e-mail: $\{$holger.dette, kevin.kokot$\}$@rub.de}\\
}

  \maketitle

\begin{abstract} 
We study the problem of testing the equivalence of functional parameters (such as the mean or variance function) in the two sample  functional  data problem.
In contrast to previous work, which reduces the functional problem to a  multiple  testing problem for the equivalence  of scalar data by comparing the functions at each point, our approach is based on an estimate
of a distance measuring the maximum deviation between  the  two functional parameters. Equivalence is claimed if the  estimate for the 
maximum deviation does not exceed a given threshold. A bootstrap procedure is proposed to obtain quantiles for the distribution of  the test statistic 
and consistency of the corresponding  test is proved in the large sample scenario.
As the methods proposed here avoid the use of the intersection-union principle they are  less conservative  and more powerful than the currently available methodology. 
\end{abstract}


\section{Introduction}  \label{sec:Intro}
\def\theequation{1.\arabic{equation}}
\setcounter{equation}{0}

Equivalence tests are nowadays frequently  used in   drug development to assess similarity of a test and  a reference treatment
at a controlled type I error. They are very popular in regulatory settings because they reverse the burden of proof compared 
to a standard test of significance. Therefore they avoid  the problem that  failing to reject a null hypothesis of no difference is not logically equivalent 
to deciding for the null hypothesis.
 Typically equivalence testing is  based on a null hypothesis that a scalar  parameter of interest, 
such as the effect difference 
of two treatments, is outside an equivalence region defined through an appropriate choice of an interval
depending on the metric of equivalence being used.  Thus rejecting the null hypothesis means 
to decide at a controlled type I error that the parameter of interest is in the postulated equivalence region.
We refer to the  monographs of \cite{wellek2010testing} 
for an overview of the currently available methodology on testing the equivalence of finite dimensional parameters.  

On the other hand there are many applications, where the similarity 
between two populations cannot be appropriately described by a  parameter of finite dimension. One obvious situation  occurs if  treatments involving covariates have to be compared and one is interested in the similarity 
of  the  relations between  the measured endpoints and  the covariates in  the two groups. Statistically speaking,
this corresponds to the problem establishing the similarity between two regression models and in the last decade considerable efforts
have been made to develop methodology to solve this problem.
 \cite{liubrehaywynn2009} proposed tests for the hypothesis of equivalence of two linear  regression models,
 while  \cite{gsteiger2011}  developed   a bootstrap  approach  using a   confidence band for the difference of two non-linear models.
  These methods are based  on 
  the intersection-union principle   \citep[see, for example,][]{berger1982} which is  used to construct an overall  test for equivalence.
In a recent paper  \cite{detmolvolbre2015}  showed    that equivalence tests  based on the intersection-union principle  
lead  to rather conservative decision procedures with low power. As a very powerful  alternative  they  proposed  bootstrap 
tests based on estimates of the  maximal deviation between the two  curves  corresponding to the different treatments.  
\cite{moedetkotvolcol2019} demonstrated the superiority of the maximum deviation
approach for the  comparison  of dissolution profiles of two different formulations  
 \citep[see][for some alternative equivalence tests based  
 on  similarity factors]{paixao2017,yoshida2017}. In all these papers, data is finite dimensional and the curves to be compared  
are defined by parametric regression models with finite dimensional parameters. 

Moreover, in the information age, data is often recorded sequentially over time at high resolution and in such instances it is reasonable
to model data as functions because  the  densely sampled observations exhibit certain degrees of dependence and smoothness.   As a consequence 
corresponding parameters  such as mean or variance  are varying over time and have to be considered as functions
 as well.  The corresponding
field in statistics is called  functional data analysis  and the current state of the art in analyzing functional data is well documented in the monographs
by  \cite{RamsaySilverman2005},  \cite{FerratyVieu2010}, \cite{HorvathKokoskza2012}, and \cite{hsingeubank2015}. Although numerous 
statistical concepts such as the comparison of mean functions, covariance operators, principal components,  change point analysis have been  
considered and developed for functional data, the problem of establishing the  practical equivalence of two parameters (more precisely parameter functions) 
for functional data has not found much attention in the literature. 

In a recent paper  \cite{fogarty2014} developed methodology for establishing the equivalence 
between the mean  and variance functions from two populations. Their work is motivated by a  comparison study  of  devices for assessing pulmonary function
and extends the  popular Two One-Sided Testing (TOST) procedure for equivalence testing of scalars \citep[see][among others]{schuirmann1987,phillips1990} to the 
functional regime. 
By the duality between  hypotheses testing and confidence intervals  their approach is equivalent to the construction of a lower and an upper  (pointwise) confidence band  
for the difference of the two parameters.
The test then decides for equivalence if the functions $\kappa_{l}(\cdot )$ and $\kappa_{u}(\cdot )$ defining the lower and upper  equivalence region for the difference of the two functional parameters 
 are outside  of the upper and lower confidence band. Thus their method  is similar in spirit to the work of \cite{liubrehaywynn2009} and   \cite{gsteiger2011}  for the comparison of parametric regression models and therefore 
expected  to be rather conservative.
A similar comment applies to equivalence tests that  can be  constructed in the same way using simultaneous confidence bands  as 
developed in  \cite{dette2018} and \cite{lieblreim2019}.

The purpose of this paper is to develop more efficient procedures to establish equivalence of parameters  
for the  two sample problem in functional  data analysis.  Our approach is based on an estimate of the maximum deviation between 
parameter functions (such as the difference of the mean functions  or the ratio of the variance functions) and we propose to decide for similarity if the estimated distance is small. In Section \ref{sec2} we introduce the basic model and 
review the method  of  \cite{fogarty2014}. Section~\ref{sec3} is devoted to the construction of a more powerful  test for the  equivalence of functional parameters,
where we concentrate on the mean functions for the sake of brevity. In particular a bootstrap test is developed and its consistency  is proved. We also provide a generalization to dependent data and illustrate the superiority of the new test in a small example. In Section \ref{sec4} we 
demonstrate the general applicability of our approach and  develop methodology for a functional random effect model as considered by  \cite{fogarty2014}.
We also demonstrate by means of a simulation study that the new tests introduced in this paper are more powerful than the currently available methodology. 
Finally, all proofs are given in an appendix as they are technically demanding and involve functional data analysis for Banach space valued random variables.

 \section{Formulation of the problem and state of the art} \label{sec2} 
 \def\theequation{2.\arabic{equation}}
\setcounter{equation}{0}

In this section we state the problem and briefly revisit the approach proposed in \cite{fogarty2014}. To
be precise, let $X_{11}(\cdot),\ldots,X_{1m}(\cdot)$ and $X_{21}(\cdot),\ldots,X_{2n}(\cdot)$ denote two independent samples of functional data, which are observed on the interval $[0,1]$. We denote the mean functions by  $\mu_1(\cdot)$ and $\mu_2(\cdot)$ and variance functions by $\sigma^2_1(\cdot)$ and $\sigma^2_2(\cdot)$,  respectively (assuming its existence - see Section \ref{sec6} for the necessary assumptions). 
 We define $\theta (\cdot) = \mu_1(\cdot) - \mu_2(\cdot)$ and $\lambda (\cdot) = \frac {\sigma^2_1(\cdot)}{\sigma^2_2(\cdot)}$ as measures of similarity 
between the mean and variance functions, respectively, and consider the hypotheses
\begin{align} 
\label{H0mean}
\begin{split}
  H^\theta_0  : & ~\exists \ t \in [0,1] \mbox { such that } \theta(t)  \notin (\kappa_l(t), \kappa_u(t)) \\ 
  H^\theta_1  : & ~\forall \ t \in [0,1]:  \theta(t) \in (\kappa_l(t), \kappa_u(t))
  \end{split}
\end{align}
and
\begin{align} 
\label{H0var}
\begin{split}
  H^\lambda_0  : & ~\exists \ t \in [0,1] \mbox { such that } \lambda(t)  \notin (\zeta_l(t), \zeta_u(t)) \\ 
  H^\lambda_1  : & ~ \forall \ t \in [0,1]: \lambda(t) \in (\zeta_l(t), \zeta_u(t)) \, .
\end{split}
\end{align}
Here $\kappa_l, \kappa_u, \zeta_l, \zeta_u$ are given functions  on the interval $[0,1]$, which define the region of equivalence. These bands have to be developed in cooperation with the experts from the field of application.  Usually  the band defined by the functions $\kappa_l $ and $ \kappa_u$ contains the constant function 
 $0$ (as one wants to demonstrate the similarity of the functions $\mu_1$ and $\mu_2$) and the band defined by the functions 
$\zeta_l$ and $\zeta_u$ contains the constant function $1$. 
Note that the rejection of the null hypothesis in \eqref{H0mean} means to decide (at a controlled type I error) that the difference of the mean functions is contained 
in the band defined by the functions $\kappa_l$ and $ \kappa_u$ and a similar comment applies  to the rejection of
 the null hypothesis in \eqref{H0var}.

In the following, we concentrate on the mean functions to describe the currently available methodology.
\cite{fogarty2014} combined the intersection-union principle with equivalence testing of scalar parameters to develop tests for the hypotheses \eqref{H0mean}.  More precisely, they 
 proposed to test for equivalence in location at each $ t \in [0,1]$ and to reject the null hypothesis in \eqref{H0mean} if all individual tests yield a rejection. For the construction of 
 the individual tests  they used a bootstrap version of the Two-One-Sided-Testing (TOST) principle as introduced by  \cite{schuirmann1987}. 
 To be precise, if $\hat \theta_{m,n} (t) = \overline{X}_{1 \cdot}(t) - \overline{X}_{2 \cdot}(t)$ is the common estimate of the mean difference at time $t \in [0,1]$ and
\begin{align*}
  \overline{C} _{1-\alpha,\theta} (t) &= \big[ 2 \hat\theta_{m,n}(t) 
  - q_{1-\alpha} (\hat\theta^*_{m,n}(t)), \infty \big) \\
  \underline{C}_{1-\alpha,\theta} (t) &= \big (- \infty, 2 \hat \theta_{m,n}(t) 
  - q_{ \alpha}(\hat \theta^*_{m,n}(t)) \big ]
\end{align*}
are bias corrected percentile-based bootstrap one-sided confidence intervals, then the individual null hypothesis $H_{0,t}^\theta :\theta(t) \notin (\kappa_l (t), \kappa_u(t))$ is rejected in favor of
$H_{1,t}^\theta: \theta(t) \in (\kappa_l(t), \kappa_u(t))$ if $\kappa_l(t) \notin \overline{C}_{1-\alpha,\theta}(t)$ 
{\bf and} $\kappa_u \notin \underline{C}_{1-\alpha,\theta}(t)$, or equivalently
\begin{equation}\label{hd3}
  \kappa_l(t) < 2 \hat \theta_{m,n}(t) - q_{1-\alpha}(\hat \theta^*_{m,n}(t)) 
  \leq 2 \hat \theta_{m,n}(t) - q_\alpha (\hat \theta^*_{m,n}(t)) 
  < \kappa_u(t) \, .
  \end{equation}
  
  \smallskip

\begin{rem}  \label{rem1}~
{\rm 
\begin{itemize}
\item[(a)]
It is worthwhile to mention that  the concept described here can be used with any type of one-sided confidence intervals.
 For example, it follows  from the proofs of the results in Section~\ref{sec3} (see Section \ref{sec6} for more details)  that 
$\sqrt{m+n} ( \hat \theta_{m,n}(t) -  \theta(t) ) $ (for fixed  $t$) is asymptotically normal distributed with  variance $\sigma^{2} (t) = \frac{1}{\tau }\sigma^2_1(t) + \frac {1}{1-\tau }\sigma^2_2(t)$,
 where $\tau = \lim_{m,n \to \infty} m/(m+n)$. Therefore  one could use the asymptotic $(1-\alpha)$-confidence intervals
 \begin{align*}
   \overline{C}_{1-\alpha,\theta}^{a} (\theta(t)) 
   &= \Big  [\hat \theta_{m,n}(t) - u_{1-\alpha} \frac{\hat \sigma_{m,n} (t)}{{\sqrt{m+n}}} \, , \infty \Big) \\
   \underline{C}_{1-\alpha,\theta}^{a} (\theta(t) ) 
   &= \Big (- \infty, \hat \theta_{m,n}(t)-u_{\alpha} \frac{\hat \sigma_{m,n} (t)}{{\sqrt{m+n}}} \, \Big]
\end{align*}
 to derive an analogue of the decision rule \eqref{hd3}, where $\hat \sigma_{m,n}^{2} (t)$ is an appropriate estimate of the asymptotic variance $\sigma^{2} (t) $ at the point $ t \in [0,1]$ and $u_\alpha$ denotes the $\alpha$-quantile of the standard normal distribution.
 \item[(b)] Besides the  frequentist test described in the previous  paragraph \cite{fogarty2014} also proposed a test within the Bayesian 
 paradigm using Gaussian Processes for  modelling the data. Because the focus of this paper is on nonparametric procedures we do not consider this test here.
 \end{itemize}
 }
\end{rem}

\section{Efficient equivalence-testing of functional parameters } \label{sec3} 
\def\theequation{3.\arabic{equation}}
\setcounter{equation}{0}

In this section we develop an alternative test for the hypotheses \eqref{H0mean} in the two sample problem, which turns
out to be substantially more powerful than the frequentist method proposed by \cite{fogarty2014}.
  Our approach is based 
on the estimation of the maximum deviation of the unknown measure of similarity from the equivalence bounds defined
in \eqref{H0mean} and \eqref{H0var}. 
To be precise  we restrict ourselves  again to the difference of the location parameters
$\theta = \mu_{1} - \mu_{2} $ and note that the hypotheses in \eqref{H0mean} can be rewritten as 
\begin{align} \label{eq:equi-hypotheses}
\begin{split}
  H_0^{\theta} : ~ T^{\theta}
  & = \max\Big\{  \sup_{t\in [0,1]} \big(-\theta(t) + \kappa_l(t) \big), \, 
  \sup_{t\in [0,1]} \big(\theta(t) - \kappa_u(t) \big) \Big\} \geq 0  \\
  H_1^{\theta} :~ T^{\theta } 
  &=  \max\Big\{  \sup_{t\in [0,1]} \big(-\theta(t) + \kappa_l(t) \big), \, 
  \sup_{t\in [0,1]} \big(\theta(t) - \kappa_u(t) \big) \Big\}< 0 \, .
\end{split}
\end{align}  
The representation of the hypotheses  simplifies  in the case of symmetric and constant boundaries, that is  $\kappa_u(t) =- \kappa_l (t) = \kappa > 0 $ for all $t \in [0,1]$, where we obtain
for the hypotheses in \eqref{eq:equi-hypotheses}
$$
 H_0^{\theta} :  \sup_{t\in [0,1]} \big |  \theta(t)   \big |  \geq  \kappa \, ,~~   H_1^{\theta} :  \sup_{t\in [0,1]} \big |  \theta(t)   \big |  <  \kappa \, .
$$ 
For the construction of an efficient test, we define the statistic
\begin{align} \label{eq:statistic}
  \hat{T}_{m,n}^{\theta} = \max\Big\{ 
  \sup_{t\in [0,1]} \big(-\hat{\theta}_{m,n}(t) + \kappa_l(t) \big), \,
  \sup_{t\in [0,1]} \big(\hat{\theta}_{m,n}(t) - \kappa_u(t) \big) \Big\}
\end{align}
as an estimator of $T^{\theta}$,  where 
$$
\hat{\theta}_{m,n} = \overline{X}_{1\cdot}  - \overline{X}_{2\cdot} 
$$
 denotes the difference of  the sample means, which serves as an estimator 
 of  the function  $\theta= \mu_1-\mu_2$.  The null hypothesis in \eqref{eq:equi-hypotheses} is then rejected for small values
 of $  \hat{T}_{m,n}^{\theta}$, where the critical values will be determined by  bootstrap (in the independent case by resampling with replacement, in the dependent case by multiplier block bootstrap). 
 
 To be precise and motivate our bootstrap assume that  $m,n \to \infty$, such that $m/(m+n) \to \tau \in (0,1)$. Then
if follows from Theorem \ref{thma1} in the online supplement that 
 \begin{align} \label{eq:limit-dist}
\sqrt{m+n} \big ( \hat T^\theta_{m,n} - T^\theta \big) {\stackrel{\mathcal{D}}{\longrightarrow}}
  Z_{\mathcal{E}, \theta} = \max\Big\{ 
  \sup_{t\in \mathcal{E}^l_ \theta} \big(-Z(t)\big), \,
  \sup_{t\in \mathcal{E}^u_ \theta} Z(t) \Big\} \, ,
\end{align}
where $Z$ is a Gaussian process with covariance kernel  
\begin{align} \label{eq:Z-kernel}
k(s,t) =  \frac{1}{\tau} \mbox{Cov} (X_{11}(s), X_{11}(t)) + \frac{1}{1-\tau} \mbox{Cov} (X_{21}(s), X_{21}(t))
\end{align}
and the  sets $ \mathcal{E}^l_{ \theta} ,  \mathcal{E}^u_ {\theta} \subset [0,1] $ contain the points, where 
the functions  $- \theta +\kappa_l $  and $\theta -\kappa_u $ attain the  value $T^{\theta } $,  i.e. 
  \begin{align} \label{eq:sets}
  \mathcal{E}^l_\theta  = 
  \big\{ t\in [0,1] \colon -\theta(t)+\kappa_l(t) = T^{\theta } \, \big\} \, , \quad
  \mathcal{E}^u_\theta   = 
  \big\{ t\in [0,1] \colon \theta(t)-\kappa _u(t) = T^{\theta } \, \big\} \, . 
\end{align}
Throughout this paper, these sets are called {\it  extremal sets} and we 
note that the extremal sets can be empty (but not both at the same time). 
As a consequence, the limit distribution on the right hand side of \eqref{eq:limit-dist} depends on the covariance kernel $k$  and 
the extremal sets $  \mathcal{E}^l_\theta$ and $  \mathcal{E}^u_\theta$ defined by the 
unknown difference $\theta$ between the mean functions $\mu_{1}$ and $\mu_{2}$.

For the calculation of quantiles  of  the distribution of $ Z_{\mathcal{E}, \theta}$ we propose to use the bootstrap and  proceed in two steps:
\begin{itemize}
\item[(1)]  We estimate the unknown sets of extremal points.
\item[(2)]  We use the bootstrap to mimic the distribution of the process $Z$ in \eqref{eq:limit-dist}. 
\end{itemize}
For the estimation of the extremal sets $  \mathcal{E}^l_\theta $  and $  \mathcal{E}^u_\theta $, we use the statistics 
\begin{align} \label{eq:est-sets0}
\begin{split}
  \hat{\mathcal{E}}_{\theta}^l  &= 
  \Big\{ t\in[0,1] \colon -\hat{\theta}_{m,n}(t)+\kappa_l(t)
  \geq \hat{T}_{m,n}^\theta - c \, \frac{\log(m+n)}{\sqrt{m+n}} \, \Big\} \, , \\
  \hat{\mathcal{E}}_{\theta}^u  &= 
  \Big\{ t\in[0,1] \colon \hat{\theta}_{m,n}(t)-\kappa_u(t)
  \geq \hat{T}_{m,n}^\theta - c \, \frac{\log(m+n)}{\sqrt{m+n}} \, \Big\} \, ,
\end{split}
\end{align}
where the statistic $ \hat{T}_{m,n}^\theta$ is defined in \eqref{eq:statistic} and $c$ is a tuning parameter.  
 For the bootstrap part  note that it follows from the  arguments  given 
 in the proof of Theorem~\ref{thma1} in the appendix that the statistic on the  left hand side of   \eqref{eq:limit-dist} 
is asymptotically equivalent to the statistic
$$
\max\Big\{ 
  \sup_{t\in \mathcal{E}^l_ \theta} \big(-\hat Z_{m,n}  (t)\big), \,
  \sup_{t\in \mathcal{E}^u_ \theta}   \hat Z_{m,n} (t) \Big\}
$$
where the  process $\hat Z_{m,n} $ is defined by 
$$
 \hat Z_{m,n}= \sqrt{m+n}  \, \big\{ \hat \theta_{m,n} - \theta \big\}
 = \sqrt{m+n}  \, \big\{ \overline{X}_{1 \cdot} - \mu_{1} 
 - ( \overline{X}_{2 \cdot} -\mu_{2} ) \big\} 
$$
(by the arguments given in the proof of Theorem~\ref{thma1}  this process converges weakly to the process $Z$ on the right 
hand side of \eqref{eq:limit-dist}).
To mimic the distribution of this process in the independent case we now  use   resampling with replacement.
More precisely,  assume  for $r=1, \ldots , R$ that 
$X^{*(r)} _{11}, \ldots, X^{*(r)}_{1m}$ and 
$X^{*(r)}_{21}, \ldots, X^{*(r)}_{2n}$ are drawn randomly with replacement
from $X_{11},\ldots, X_{1m}$ and $X_{21},\ldots,X_{2n}$, respectively, and 
denote by $\overline{X}_{1 \cdot}$ and $\overline{X}_{2 \cdot}$ the sample 
means of both groups. We define 
\begin{align} \label{eq:bootstrap-process}
  \hat Z_{m,n}^{*(r)} =
  \sqrt{m+n} \, \Big \{ \frac {1}{m} \sum^m_{j=1} 
  (X^{*(r)}_{1j} - \overline{X}_{1 \cdot}) 
  - \frac {1}{n} \sum^n_{j=1} (X^{*(r)}_{2j}- \overline{X}_{2 \cdot}) \Big \}
\end{align}
as the $r$-th bootstrap analogue of the statistic $ \hat Z_{m,n}$
and a  bootstrap version of the random variable on the left  hand side of \eqref{eq:limit-dist} by  
\begin{equation} \label{hda}
\hat T^{\theta, * (r)}_{m,n} = \max \Big \{ \sup_{t \in \hat {\mathcal{E}}^l_\theta} \big ( - \hat Z^{*(r)}_{m,n}(t) \big), \, \sup_{t \in \hat{\mathcal{E}}^u_\theta} \hat Z^{*(r)}_{m,n} (t) \Big \} \, .
\end{equation}
Finally, the null hypothesis in \eqref{eq:equi-hypotheses}  is  rejected, whenever
\begin{align} \label{eq:equi-test}
 \sqrt{m+n} \, \hat{T}_{m,n}^\theta < z_{m,n,\alpha}^{*(R)} \, ,
\end{align}
where $z_{m,n,\alpha}^{*(R)}$ is the empirical $\alpha$-quantile of the 
bootstrap sample $ T^{\theta,*(1)}_{m,n}, \ldots, T^{\theta,*(R)}_{m,n}$. 
The following result, which is proved  in the appendix, shows  that this procedure defines a consistent and 
asymptotic level $\alpha$-test for the hypotheses  \eqref{eq:equi-hypotheses} 
(or equivalently for the hypotheses \eqref{H0mean}).

\begin{theorem} \label{thm1}
Let  Assumption~\ref{as:ts} in Section~\ref{sec62} be satisfied.
\begin{itemize}
\item[(a)] Assume that the null hypothesis $H_{0}^{\theta}$ of no 
equivalence in \eqref{H0mean} holds, that is 
$ T^{\theta} \geq 0$. 
 If $T^\theta =0$, then  
 \begin{align*}  
  \lim_{m,n,R\to\infty} \mathbb{P}\big( 
  \sqrt{m+n} \, \hat{T}_{m,n}^\theta < z_{m,n,\alpha}^{*(R)} \big) = \alpha \, .
\end{align*}
 If $T^\theta > 0$, then for any $R \in \N$
\begin{align*} 
 \lim_{m,n\to\infty} \mathbb{P}\big( 
  \sqrt{m+n} \, \hat{T}_{m,n}^\theta < z_{m,n,\alpha}^{*(R)} \big) = 
    0 \, .
\end{align*}
\item[(b)] If the alternative  $H_{1}^{\theta}$  of equivalence in \eqref{H0mean}  holds, that is  $T^\theta < 0$, we have for any $R \in \N$
\begin{align*} 
\liminf_{m,n\to\infty} \mathbb{P}\big ( \sqrt{m+n} \, \hat{T}_{m,n}^\theta < z_{m,n,\alpha}^{*(R)} \big)
=1 \, .
\end{align*}
\end{itemize}
\end{theorem}

\begin{rem} \label{rem3} 
{\rm 
The results remain correct in the case of dependent data, where $X_{1,1}, \ldots  X_{1,m}$ and $X_{2,1}, \ldots  X_{2,n}$ are  two independent stationary time series.
In this case, we propose to use a block multiplier bootstrap to mimic the dependency in the data. To be precise, 
define a  bootstrap process  by
\begin{align} \label{2bProcess}
\begin{split}
\hat Z_{m,n}^{**(r)}(t) =& \sqrt{m+n} \Big\{
\frac{1}{m} \sum_{k=1}^{m-l_1+1} \frac{1}{\sqrt{l_1}}\Big( \sum_{j=k}^{k+l_1-1} X_{1j}(t)
-\frac{l_1}{m}\sum_{j=1}^m X_{1j}(t) \Big) \xi_k^{(r)} \\
&- \frac{1}{n} \sum_{k=1}^{n-l_2+1} \frac{1}{\sqrt{l_2}}\Big( \sum_{j=k}^{k+l_2-1} X_{2j}(t)
-\frac{l_2}{n}\sum_{j=1}^n X_{2j}(t) \Big) \zeta_k^{(r)} \Big\}  ~~~~(r=1, \ldots , R)
\end{split}
\end{align}
where $\xi_1^{(r)} , \ldots , \xi_m^{(r)}$, $ \zeta_1^{(r)} , \ldots , \zeta_n^{(r)} $ are independent standard normal distributed random variables
and $l_{1}, l_{2}$ are sequences converging to infinity
with increasing sample sizes $m,n \to \infty$.
The null hypothesis  in \eqref{H0mean} is  now rejected, whenever
\begin{align} \label{eq:equi-testdep}
 \sqrt{m+n} \, \hat{T}_{m,n}^\theta < z_{m,n,\alpha}^{**(R)}~,
\end{align}
where $z_{m,n,\alpha}^{**(R)}$ is the  empirical $\alpha$-quantile  of the sample  $  T^{\theta,**(1)}_{m,n}, \ldots, T^{\theta,**(R)}_{m,n}$ 
and the statistic $  \hat T^{\theta, ** (r)}_{m,n}$ is defined by
\begin{align*}
  \hat T^{\theta, ** (r)}_{m,n} = \max \Big \{ \sup_{t \in \hat {\mathcal{E}}^l_\theta} \big ( - \hat Z^{**(r)}_{m,n}(t) \big), \, \sup_{t \in \hat{\mathcal{E}}^u_\theta} \hat Z^{**(r)}_{m,n} (t) \Big \}  ~~~~(r=1, \ldots , R) \, .
\end{align*}
In this case - under the assumptions stated in Section \ref{sec623} -  the result in Theorem \ref{thm1} remains valid.
Finally we note that this procedure with $l_{1} = l_{2}=1$ provides also a valid bootstrap test in the case of independent data.

}
\end{rem}

\begin{example} \label{examsim1} 
{\rm 
We have conducted  a small 
simulation study  to compare the new bootstrap test \eqref{eq:equi-test} with  the 
frequentist test proposed by \cite{fogarty2014}.  Further numerical results supporting our findings can be found 
in Section  \ref{sec53}, where we compare both methods in a functional random effect model.

We have  generated  functional data  as described in Sections~6.3 and 
6.4 of \cite{aueDubartNorinhoHormann2015}, who considered $D =21 $ 
$B$-spline basis functions $\nu_1,\dots , \nu_D$ and defined the 
random functions $\eta_{11}, \dots, \eta_{1m}, \eta_{21}, \dots, \eta_{2n}$ by 
\begin{equation}\label{eq:errors}
  \eta_{1j} = \sum_{i=1}^D N_{1,i,j} \nu_i \, , \quad 
  \eta_{2k} = \sum_{i=1}^D N_{2,i,k} \nu_i \, ,
  \qquad j=1,\dots, m, \, k=1,\dots,n \, ,
\end{equation}
where $N_{1,1,1}, \dots, N_{1,D,m}, N_{2,1,1}, \dots, N_{2,D,n} $ are 
independent, normally distributed random variables with expectation zero and 
variances $\sigma_i^2 = \var(N_{1,i,j}) = \var(N_{2,i,k}) = 1/i^2$ 
($i = 1,\dots, D$; $j=1,\ldots , m$; $k=1,\ldots , n$).

\begin{figure}[t]
  {  \centering
    \includegraphics[width = 6cm, height = 6cm]{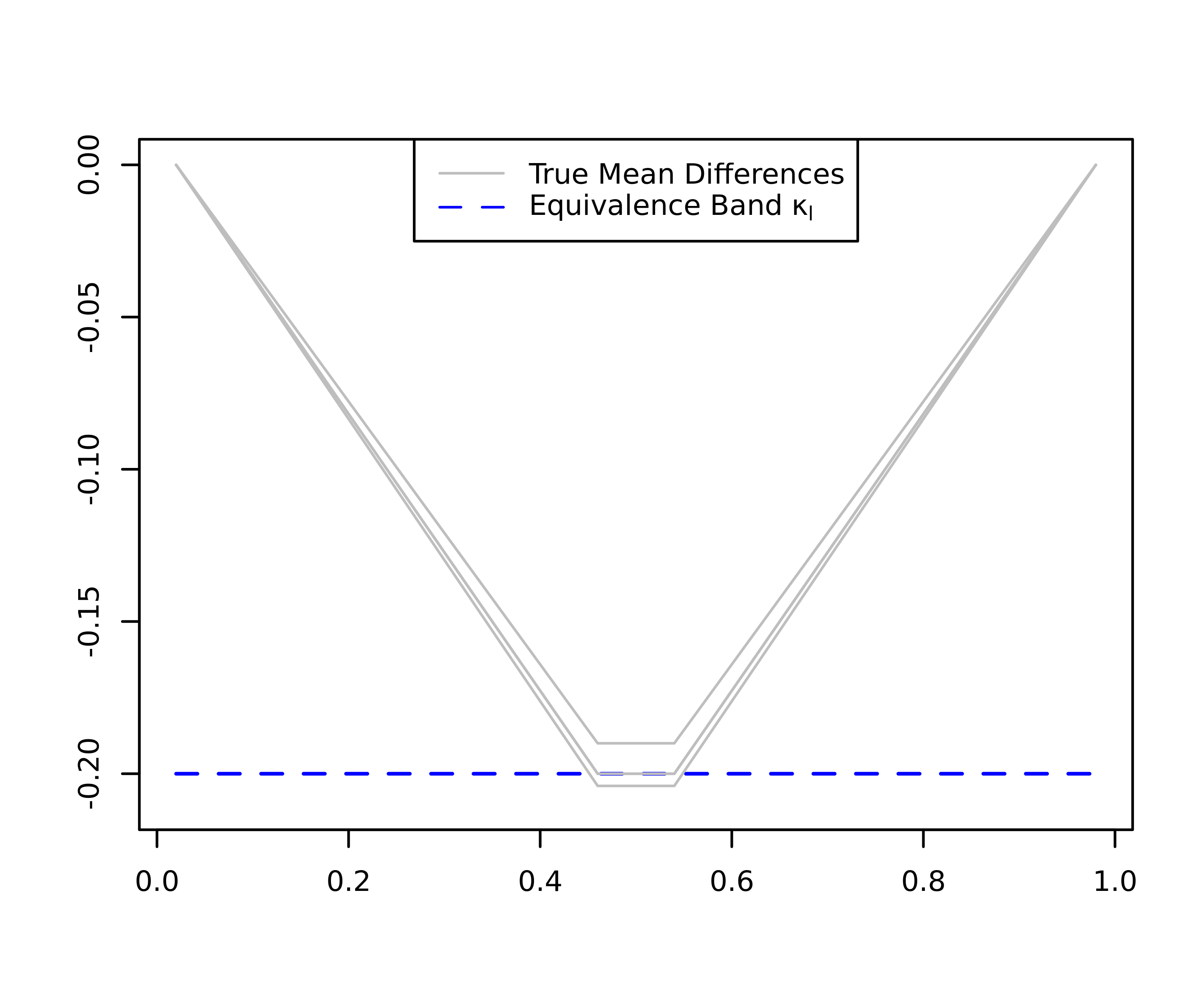}
    \includegraphics[width = 6cm, height = 6cm]{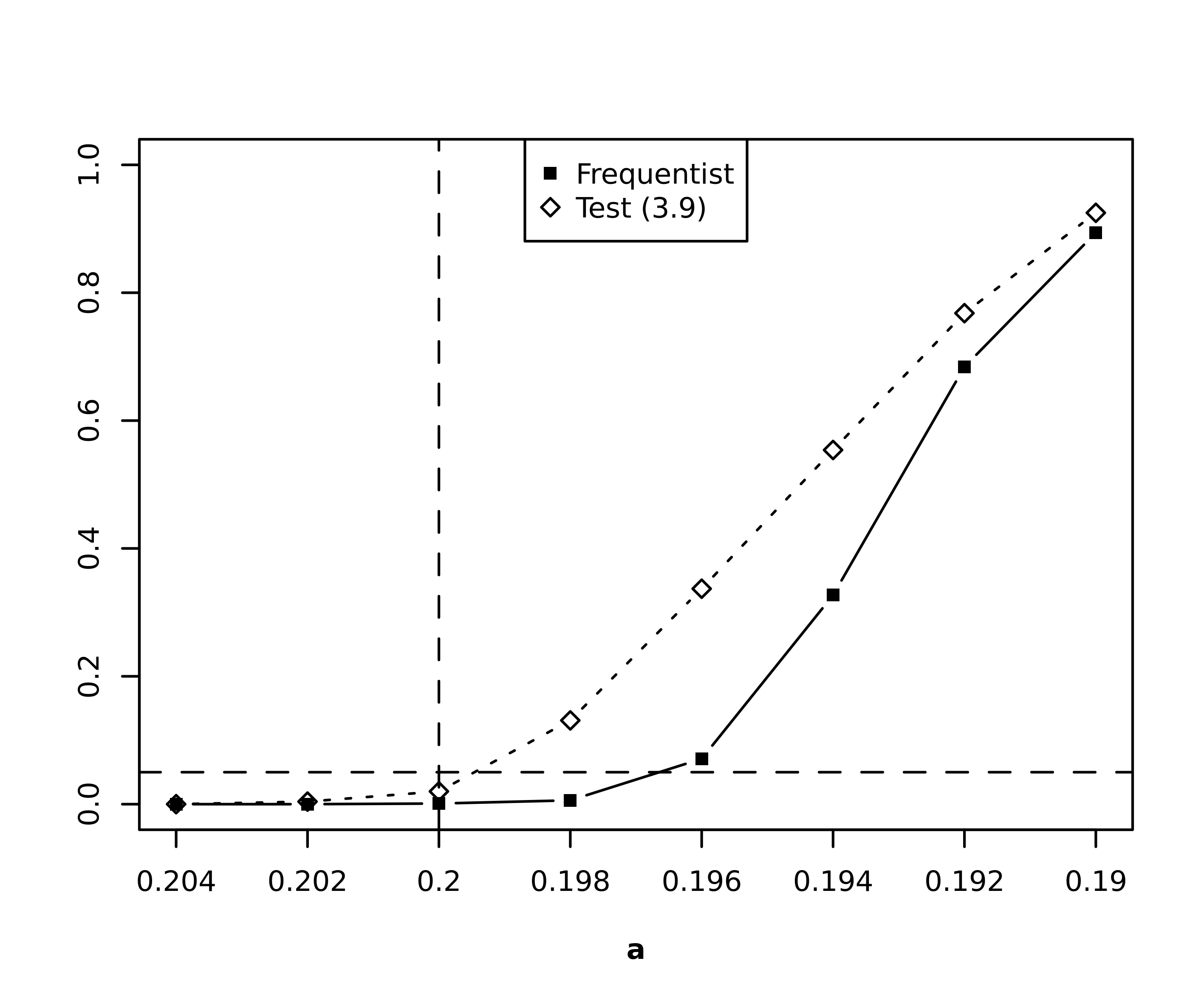}
    \caption{\it Left panel: Difference $\theta = \mu_1 -  \mu_2$ of the mean 
      functions defined by \eqref{eq:subinterval} with fixed $b_1 = 0.46, b_2 = 0.54$ 
      and different values for $a\in \{ 0.19,  0.204, 0.2 \}$.
      Right panel: Empirical rejection probabilities of the frequentist test  
      proposed by \cite{fogarty2014} and the test \eqref{eq:equi-test} for the  
      hypotheses \eqref{eq:equi-hypotheses}  with 
      $\kappa_{l}  \equiv -0.2$, $\kappa_{u} \equiv 0.2$ and different values of $a$.
    }
    \label{fig:subinterval}}
\end{figure}

Then two independent samples of  independent and identically distributed Gaussian random 
functions are obtained by 
\begin{align} \label{eq:fIID}
X_{1j} = \eta_{1,j} + \mu_1 \, , \quad X_{2k} = \eta_{2,k} + \mu_2 \, ,
\end{align} 
with  mean functions 
\begin{align} \label{eq:subinterval}
\mu_1 \equiv 0 \, , \quad \mu_2(t) = 
\begin{cases}
\frac{a}{b_1-0.02} (t-0.02) \, , & t\in [0, b_1) \\
a \, ,        & t\in[b_1, b_2] \\
\frac{-a}{0.98 - b_2} (t - b_2) + a \, , &t \in (b_2,1]
\end{cases}~,
\end{align}
where    $a$,  $b_1$ and $b_2 $ are parameters.
The left part of Figure~\ref{fig:subinterval} illustrates  the difference $\theta = \mu_{1} - \mu_2$ of the mean functions  for fixed $b_1 = 0.46$, $b_2 = 0.54$ and different values of the parameter $a$.
The equivalence bands, used in the hypotheses \eqref{H0mean}, are defined by  
$\kappa_l \equiv -0.2$, $\kappa_u \equiv 0.2$. Note that, for any $a > 0.02$, 
the extremal sets in \eqref{eq:sets} are defined by $\mathcal{E}_\theta^u = \emptyset$,
$\mathcal{E}_\theta^l = [b_1, b_2]$ (here 
$\mathcal{E}_\theta^l = [0.46,0.54]$) and that the cases $|a| \geq 0.2$ and $|a| <  0.2$ correspond to the null hypothesis of no 
equivalence and   the alternative  of equivalent mean functions, respectively. In  the  right part of 
Figure~\ref{fig:subinterval} we display the  empirical rejection probabilities 
of the frequentist test proposed by \cite{fogarty2014} and the test defined 
in \eqref{eq:equi-test}  for different values of $a \in \{  
0.204,0.202,\dots,0.190 \} $  (by symmetry negative values of $a$ yield the same 
results). Here  the extremal sets are estimated by   
\eqref{eq:est-sets0} with $c=0.005$.

The sample sizes are $m = n = 100$ and the rejection probabilities are 
calculated by $1000$ simulation runs and $300$ bootstrap replications. We 
observe that the rejection probabilities are strictly smaller than the level 
$5\%$ for $a>0.2$ and increase towards $1$ for decreasing $a$ beyond $0.2$. 
Both tests slightly underestimate the nominal level at the boundary  of the 
null hypothesis ($a=0.2$). Moreover, the new test has substantially more power 
in all considered scenarios under the alternative.

\begin{figure}[h]
  {  \centering
    \includegraphics[width = 6cm, height = 6 cm]{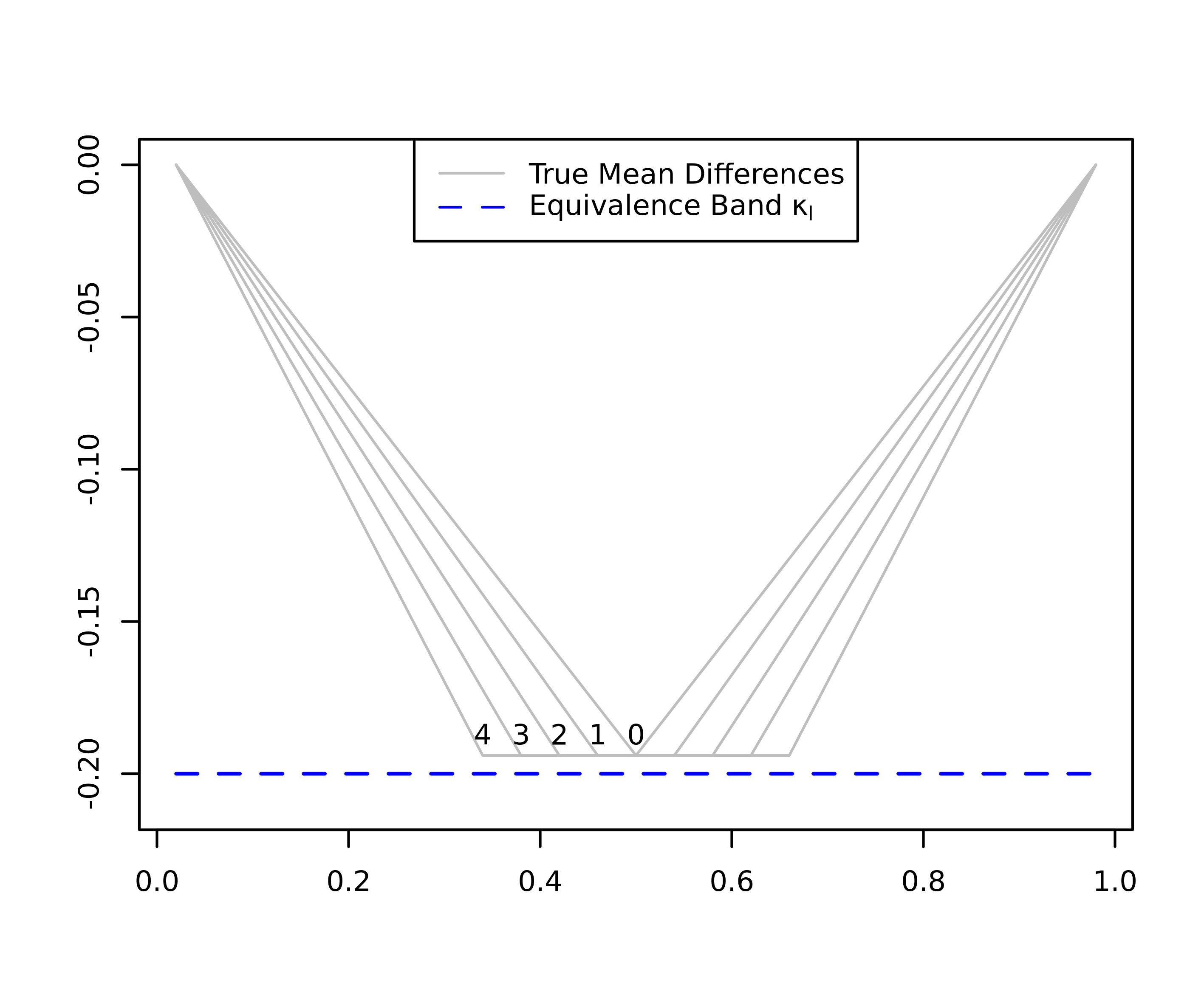}
    \includegraphics[width = 6cm, height = 6 cm]{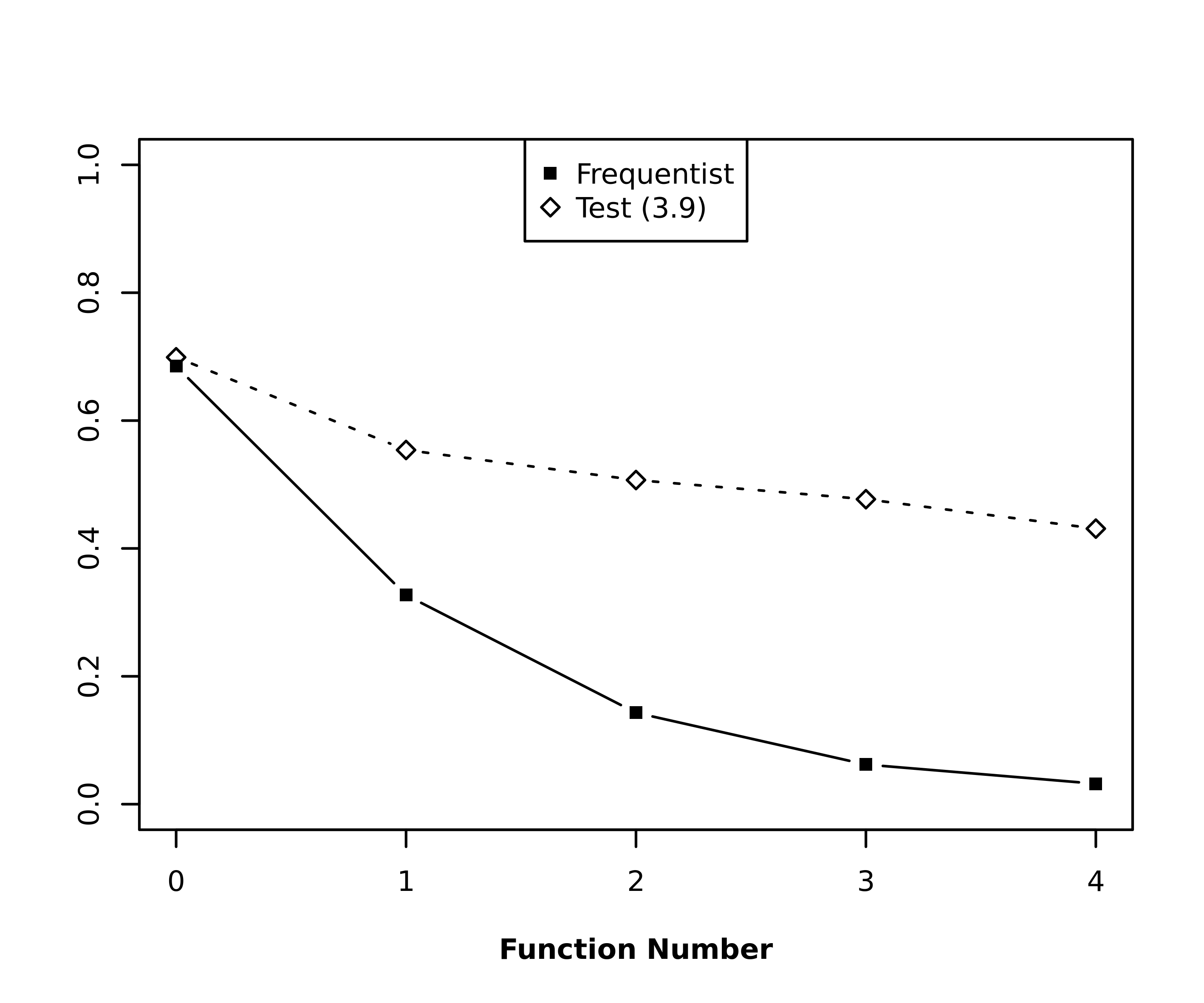}
    \caption{\it Left panel: Difference $\theta = \mu_1 - \mu_2$ of the mean 
      functions defined by \eqref{eq:subinterval} for fixed $a = 0.194$ and 
      different choices of $b_1 = 0.5-0.08j$, $b_2 = 0.5+0.08j$, where   $j=0,\dots,4$.
            Right panel: Empirical rejection probabilities of the frequentist test  proposed by \cite{fogarty2014} 
      and the test \eqref{eq:equi-test} for the  
      hypotheses \eqref{eq:equi-hypotheses}  with 
      $\kappa_{l} \equiv -0.2$, $\kappa_{u} \equiv 0.2$.
   }
    \label{fig:varySets}}
\end{figure}

The superiority of the new test is even more visible if the size of the set of 
extremal points is larger. To illustrate this fact, we consider the mean 
functions in \eqref{eq:subinterval} for fixed $a = 0.194$ and different values 
of $b_1$ and $b_2$.
The rejection probabilities of the frequentist test proposed by 
\cite{fogarty2014} and the test defined in \eqref{eq:equi-test}
are shown in Figure \ref{fig:varySets} where, for $j=0,\dots,4$, function 
number $j$ corresponds to the choices $b_1 = 0.5-0.08j$, $b_2 = 0.5+0.08j$ in 
the definition of the mean differences in \eqref{eq:subinterval}. The sample 
sizes are again $m=n=100$. We observe that only in the case $j=0$, both tests 
have comparable power. In all other cases, the new test \eqref{eq:equi-test} 
outperforms the test proposed by \cite{fogarty2014} substantially.

}
\end{example}

\section{A functional random effect model for paired data } \label{sec4}
\def\theequation{4.\arabic{equation}}
\setcounter{equation}{0}

In this section we demonstrate that the method introduced in Section \ref{sec3} for the simple two sample problem of comparing 
two mean functions is a universally applicable  decision rule to decide for  the equivalence between two functional parameters  from  two  samples 
of functional data. For this purpose only the bootstrap  procedure has to be adjusted to the situation under consideration.

As a concrete  example (in particular for the sake of comparison with the currently available methodology) we consider a functional analysis of variance model with 
random effects as  proposed by \cite{fogarty2014} for the analysis of 
functional data describing the lung volume over time for different patients and different breaths produced by a spirometer (industry standard)  and a new  device (Structured Light Plethysmography - SLP). While the new SLP holds many advantages, it has to be assured that it produces measurements (practically) equivalent to those produced by the industry standard, before it can be used for diagnoses purposes   \citep[see][for more details]{fogarty2014}. There are $A$ patients and for the $i$-th patient, $n_i$ breaths are recorded simultaneously by both devices leading to paired functional data with cross-covariances between the pairs. The goal is the development of  a statistically justified decision rule to decide  for or against equivalence of the measurements.

To be precise, we consider pairs of random functions 
defined by
\begin{align} \label{eq:fogarty-datamodel}
  \left(
  \begin{array}{ll}
    X_{1,i,j} \\
    X_{2,i,j} 
  \end{array}
  \right)
  = \left(
  \begin{array}{ll}
    \mu_1 + \varepsilon_{1,i} + \eta_{1,i,j} \\ 
    \mu_2 + \varepsilon_{2,i} + \eta_{2,i,j}
  \end{array}
  \right) \, \quad 
  j = 1,\dots,n_i \, , \, i = 1,\dots,A \, .
\end{align}
Here   $\mu_1,\mu_2$ denote the mean functions, the functions  $\varepsilon_{1,i}, \varepsilon_{2,i}$ model a random group effect
(usually corresponding to different individuals drawn from a larger population) and  the functions  $\eta_{1,i,j},\eta_{2,i,j}$ are individual random effects. 
The random group effects and the individual random effect functions are assumed to be centred and independent and identically distributed, 
respectively. Furthermore the group effects are independent of the individual ones.
Note that  the total number of pairs is given  by $N = \sum_{i = 1}^A n_i$.

\subsection{Comparing mean functions}  \label{sec51} 
For the construction of a  test for the hypotheses \eqref{eq:equi-hypotheses}, we consider the statistic 
\begin{align} \label{eq:statistic-groups}
  \hat{T}_{N}^{\theta}= 
  \sqrt{A} \, \max\big\{ 
  \sup_{t\in [0,1]} \big(-\hat{\theta}_{N}(t) + \kappa_l(t) \big), \,  
  \sup_{t\in [0,1]} \big(\hat{\theta}_{N}(t) - \kappa_u(t) \big) \big\} \, ,
\end{align}
where $\hat{\theta}_{N} =   \overline{X}_{1\cdot \cdot }  -  \overline{X}_{2\cdot \cdot } $ and 
\begin{align*}
  \overline{X}_{\ell\cdot \cdot } = \frac{1}{A} \sum_{i=1}^A \frac{1}{n_i} \sum_{j=1}^{n_i} X_{\ell,i,j}  \, ,~~\ell =1,2
\end{align*}
denote the two  sample means. The bootstrap analogue of  
\eqref{eq:statistic-groups} is defined as follows. We use the sample means 
\begin{align} \label{eq:group-mean}
  \overline{X}_{\ell, i, \cdot } 
  = \frac{1}{n_i} \sum_{j=1}^{n_i} X_{\ell,i,j}  \, ,~~
  \ell =1,2 \, ,~ i = 1,\dots,A
\end{align}
in the different groups
to estimate the  group effects by 
\begin{align*}
\hat{\varepsilon}_{1,i} = \overline{X}_{1,i,\cdot} - \overline{X}_{1,\cdot\cdot} \, , \quad 
\hat{\varepsilon}_{2,i} = \overline{X}_{2,i,\cdot} - \overline{X}_{2,\cdot\cdot} \, , ~~~(i = 1,\dots,A)  .
\end{align*}
For the bootstrap we  draw, for $r=1,\ldots , R$, samples 
$(\hat{\varepsilon}_{1,1}^{\star (r)}, \hat{\varepsilon}_{2,1}^{\star (r)}), \ldots, (\hat{\varepsilon}_{1,A}^{\star (r)}, 
\hat{\varepsilon}_{2,A}^{\star (r)})$
randomly with replacement from  the pairs 
$(\hat{\varepsilon}_{1,1}, \hat{\varepsilon}_{2,1}), \ldots, 
(\hat{\varepsilon}_{1,A},\hat{\varepsilon}_{2,A})$.
The bootstrap statistic is then defined by 
\begin{align} \label{hd21}
  \hat T_{N}^{\theta,\star  (r) } =  \max\Big\{ 
  \sup_{t\in \hat{\mathcal{E}}^l_\theta} \big(-B_N^{\star (r) }(t)\big), \,
  \sup_{t\in \hat{\mathcal{E}}^u_\theta} B_N^{\star (r) }(t) \Big\} \, ,
\end{align}
where 
\begin{align} \label{eq:bootstrap-process2}
  B_N^{\star  (r) } = \frac{1}{\sqrt{A}} \sum_{i = 1}^A 
 \big( \hat{\varepsilon}_{1,i}^{\star (r)} 
 - \hat{\varepsilon}_{2,i}^{\star (r)} \big ) \, ,
\end{align}
and the sets $\hat{\mathcal{E}}^l, \hat{\mathcal{E}}^u$ are given by
\begin{align} \label{eq:est-sets}
\begin{split}
\hat{\mathcal{E}}_{\theta}^l  &= 
\Big\{ t\in[0,1] \colon -\hat{\theta}_N(t)+\kappa_l(t)
\geq \hat{T}_N^{\theta} - c \, \frac{\log(A)}{\sqrt{A}} \ \Big\} \\
\hat{\mathcal{E}}_{\theta}^u  &= 
\Big\{ t\in[0,1] \colon \hat{\theta}_N(t)-\kappa_u(t)
\geq \hat{T}_N^{\theta}- c \, \frac{\log(A)}{\sqrt{A}} \ \Big\} \, .
\end{split}
\end{align}
The consideration of the process 
$
  B_N^{\star  (r) } $  in \eqref{eq:bootstrap-process2} is motivated by the expansion
$  \sqrt{A} \, (\hat \theta_N - \theta) 
= \frac{1}{\sqrt{A}}\sum_{i=1}^A (\varepsilon_{1,i} - \varepsilon_{2,i})
  + o_\P(1) $, which is derived in equation 
\eqref{eq:CLT2} in the online supplement.
The null hypothesis in \eqref{eq:equi-hypotheses} is finally rejected whenever 
\begin{align} \label{t1} 
\hat{T}_{N}^\theta < z_{N,\alpha}^{\star(R)}
\end{align}
where $z_{N,\alpha}^{\star(R)}$ is the empirical $\alpha$-quantile of the bootstrap  sample
$\hat T_{N}^{\theta,\star  (1) }  , \ldots , \hat T_{N}^{\theta,\star  (R) } $. The following result shows that this decision rule defines a consistent asymptotic level 
$\alpha$ test for the hypotheses in \eqref{H0mean}.

\begin{theorem} \label{thm2}
Let Assumption~\ref{(C)}  in Section~\ref{sec51-assumptions} be satisfied 
and assume that $A\to\infty$ and $\min_i^A n_i \to\infty$ as $N\to\infty$.
\begin{itemize}
  \item[(a)] Assume that the null hypothesis $H_{0}^{\theta}$ of no 
  equivalence in \eqref{H0mean} holds, that is 
  $ T^{\theta} \geq 0$. 
  If $T^\theta =0$, then  
  \begin{align*}  
  \lim_{A, \min n_i, R \to\infty} \mathbb{P}\big( 
  \hat{T}_{N}^\theta < z_{N,\alpha}^{\star(R)} \big) = 
  \alpha \, .
  \end{align*}
  If $T^\theta > 0$, then for any $R \in \N$
  \begin{align*} 
  \lim_{A, \min n_i \to\infty} \mathbb{P}\big( 
  \hat{T}_{N}^\theta < z_{N,\alpha}^{\star(R)} \big) = 
  0 \, .
  \end{align*}
  \item[(b)] If the alternative  $H_{1}^{\theta}$  of equivalence in 
  \eqref{H0mean}  holds, that is  $T^\theta < 0$, we have for any $R \in \N$
  \begin{align*} 
  \liminf_{A, \min n_i \to\infty} \mathbb{P}\big (\hat{T}_{N}^\theta < z_{N,\alpha}^{\star(R)} \big)
  =1 \, .
  \end{align*}
\end{itemize}
\end{theorem}

\subsection{Comparing variance functions}  \label{sec52} 

Recall the  definition of model \eqref{eq:fogarty-datamodel} and define 
(assuming its existence - see Section~\ref{sec61} for more details)
\begin{align*}
  \sigma_1^2(\cdot) = \E \big[\eta_{1,1,1}(\cdot)^2\big] \, , \quad 
  \sigma_2^2(\cdot) = \E \big[\eta_{2,1,1}(\cdot)^2\big] \in C([0,1])
\end{align*}
as the variance functions of the individual errors 
$\eta_{1,1,1}, \eta_{2,1,1}$. We are interested in testing the hypotheses
\eqref{H0var}, which can be rewritten as
\begin{align} \label{eq:equi-hypotheses-var}
\begin{split}
&H_0^\lambda: ~ T^\lambda
= \max\Big\{ 
\sup_{t\in [0,1]} \big(-\log \lambda(t) + \log \zeta_l(t) \big), \,
\sup_{t\in [0,1]} \big(\log\lambda(t) - \log\zeta_u(t) \big) \Big\} \geq 0 \\  
&H_1^\lambda:~ T^\lambda < 0 \, ,
\end{split}
\end{align}
where $\lambda = \frac{\sigma_1^2}{\sigma_2^2}$ is the ratio of the two variance functions and $\zeta_l(t), \zeta_u(t)$ are the  given equivalence bands. 
{Note that we} work with the logarithm of $\lambda$ to obtain  stabilized variances.  We define  
\begin{align} \label{eq:est-var}
  \hat{\sigma}_\ell^2 = \frac{1}{N-A} \sum_{i = 1}^A 
  \sum_{j=1}^{n_i} \Big (X_{\ell,i,j}
  - \frac{1}{n_i} \sum_{k=1}^{n_i} X_{\ell,i,k} \Big )^2 \, , \quad
\ell=1,2 , 
\end{align}
 estimate the variance ratio by $\hat{\lambda} = \frac{\hat{\sigma}_1^2}{\hat{\sigma}_2^2}  $ and consider  the test 
statistic 
\begin{align} \label{eq:statistic-groups_var}
  \hat{T}^{\lambda}_{N}= \sqrt{N} \, \max\Big\{  
  \sup_{t\in [0,1]} \big(- \log \hat{\lambda}(t) +  \log \zeta_l(t) \big), \, 
  \sup_{t\in [0,1]} \big(  \log  \hat{\lambda}(t) -  \log \zeta_u(t) \big) 
  \Big\} \, .
\end{align}
For the calculation of bootstrap quantiles we adapt resampling  with replacement to the random effect model \eqref{eq:fogarty-datamodel}
and estimate the individual random effects by
\begin{align*}
  \hat{\eta}_{1,i,j} = X_{1,i,j}-\overline{X}_{1,i ,\cdot} \, , \quad 
  \hat{\eta}_{2,i,j} = X_{2,i,j}-\overline{X}_{2,i ,\cdot} 
\end{align*}
for $i = 1,\dots, A$ and $j = 1,\dots, n_i$, where the group means $\overline{X}_{\ell ,i ,\cdot}$ ($i = 1,\dots, A$) are defined by \eqref{eq:group-mean}.
We now  draw with replacement $N=\sum_{i=1}^A {n_{i}}$  pairs
$(\hat{\eta}_{1,1,1}^{\star (r)}, \hat{\eta}_{2,1,1}^{\star (r)}),\dots, (\hat{\eta}_{1,A,n_A}^{\star (r)}, \hat{\eta}_{2,A,n_A}^{\star (r)})$
from  $(\hat{\eta}_{1,1,1}, \hat{\eta}_{2,1,1}),\dots, (\hat{\eta}_{1,A,n_A}, \hat{\eta}_{2,A,n_A})$ 
and  define for $r = 1,\dots,R$
\begin{align} \label{eq:bootstrap-statistic-var}
\hat T_{N}^{\lambda, \star (r)} =  \sqrt{N} \max\Big\{ 
\sup_{t\in \hat{\mathcal{E}}^l_\lambda} \big(-C_{N}^{ \star (r)}(t)\big), \,
\sup_{t\in \hat{\mathcal{E}}^u_\lambda} C_{N}^{ \star (r)}(t) \Big\}~
\end{align}
as the bootstrap analogue of \eqref{eq:statistic-groups_var},
where 
\begin{align} \label{hol11}
  C_{N}^{ \star(r)} = \frac{C_{1,N}^{ \star (r)}}{\hat \sigma_1^2} 
  - \frac{C_{2,N}^{ \star (r)}}{\hat \sigma_2^2} \, , \quad
  C_{\ell,N}^{ \star (r)} = \frac{1}{N-A} \sum_{i = 1}^A 
  \sum_{j=1}^{n_i} \big\{ (\hat{\eta}_{\ell,i,j}^{\star (r)})^2 
- \hat{\sigma}^2_\ell \big \}  \, , 
  \quad \ell=1,2
\end{align}
and
\begin{align} \label{eq:var-est-sets}
\begin{split}
\hat{\mathcal{E}}_{\lambda}^l  &= 
\Big\{ t\in[0,1] \colon -\log \hat{\lambda}(t)+ \log \zeta_l(t)
\geq \hat{T}_{N}^{\lambda} - c \, \frac{\log(N)}{\sqrt{N}} \ \Big\} \\
\hat{\mathcal{E}}_{\lambda}^u  &= 
\Big\{ t\in[0,1] \colon \log \hat{\lambda}(t) - \log \zeta_u(t)
\geq \hat{T}_{N}^{\lambda} - c \, \frac{\log(N)}{\sqrt{N}} \ \Big\} \, .
\end{split}
\end{align}
The  consideration  of the process $C_{N}^{ \star(r)}$ in \eqref{hol11} is motivated by the expansion
\begin{align*}
  \sqrt{N}(\log \hat\lambda - \log \lambda) 
  = \frac{\sqrt{N}}{N-A} \sum_{i = 1}^A 
  \sum_{j=1}^{n_i} \Big\{ \frac{(\eta_{1,i,j})^2 - \sigma^2_1}{\sigma^2_1} 
  - \frac{(\eta_{2,i,j})^2 - \sigma^2_2}{\sigma^2_2}  \Big \} + o_\P(1) \, ,
\end{align*}
which is derived in equation \eqref{eq:delta-method} in the online supplement.

Finally, the null hypothesis in \eqref{H0var}  of no equivalence  is rejected, whenever 
\begin{align} \label{eq:var-test}
\hat{T}_{N}^{\lambda} <  u_{N,\alpha}^{\star(R)}
\end{align}
where $ u_{N,\alpha}^{\star(R)}$ is the empirical $\alpha$-quantile of the sample $\hat T_{N}^{\lambda,  \star  (1)} , \ldots , \hat T_{N}^{\lambda,  \star  (R)} $.
\begin{theorem} \label{thm3}
Let Assumption~\ref{as:ts2} in Section~\ref{sec52-assumptions} be satisfied
and assume that $A\to\infty$ and $\min_i^A n_i \to\infty$ as $N\to\infty$.
\begin{itemize}
\item[(a)] Assume that the null hypothesis $H_{0}^{\lambda}$ of no equivalence 
in \eqref{eq:equi-hypotheses-var} holds, that is $T^{\lambda} \geq 0$. If 
$T^\lambda =0$, then  
\begin{align*}  
  \lim_{A , \min n_i, R \to\infty} \mathbb{P}\big( 
  \hat{T}_{N}^\lambda < u_{N,\alpha}^{\star(R)} \big) = \alpha \, .
\end{align*}
If $T^\lambda > 0$, then for any $R \in \N$
\begin{align*} 
  \lim_{A, \min n_i \to\infty} \mathbb{P}\big( 
  \hat{T}_{N}^\lambda < u_{N,\alpha}^{\star(R)} \big) = 0 \, .
\end{align*}
\item[(b)] If the alternative $H_{1}^{\lambda}$ of equivalence in 
\eqref{eq:equi-hypotheses-var}  holds, that is $T^\lambda < 0$, we have for any 
$R \in \N$
\begin{align*} 
  \liminf_{A, \min n_i \to\infty} \mathbb{P}\big ( 
  \hat{T}_{N}^\lambda < u_{N,\alpha}^{\star(R)} \big) =1 \, .
\end{align*}
\end{itemize}
\end{theorem}

\subsection{Some numerical results } \label{sec53} 

In this section we illustrate the finite sample properties of the new bootstrap 
procedure in the functional analysis of variance model 
\eqref{eq:fogarty-datamodel} and also provide a comparison with the method 
proposed in \cite{fogarty2014}. For this purpose we consider some of the 
scenarios described in Sections 10.1 and 10.2 of this reference.
As a general picture we will  demonstrate that the procedure  proposed in this paper is more powerful than the frequentist test method developed in  \cite{fogarty2014}. 
Note that \cite{fogarty2014} also develop a Bayesian test method but it is outperformed by the frequentist test. Therefore, the new bootstrap test is only compared with the frequentist test in the following sections. 

For the sake of comparison, we perform the frequentist test  of \cite{fogarty2014} with the same data as the new bootstrap procedure and do not use the exact results displayed in this reference. In each scenario under consideration, we perform $1000$ simulation runs and in each run, $R = 300$ bootstrap replicates are generated to calculate  the empirical $5\%$ bootstrap quantile. The extremal sets are estimated as in \eqref{eq:est-sets} and \eqref{eq:var-est-sets} with $c=0.005$, respectively.
  
\begin{figure}[t]
  {  \centering
    \includegraphics[width = 6cm, height = 6 cm]{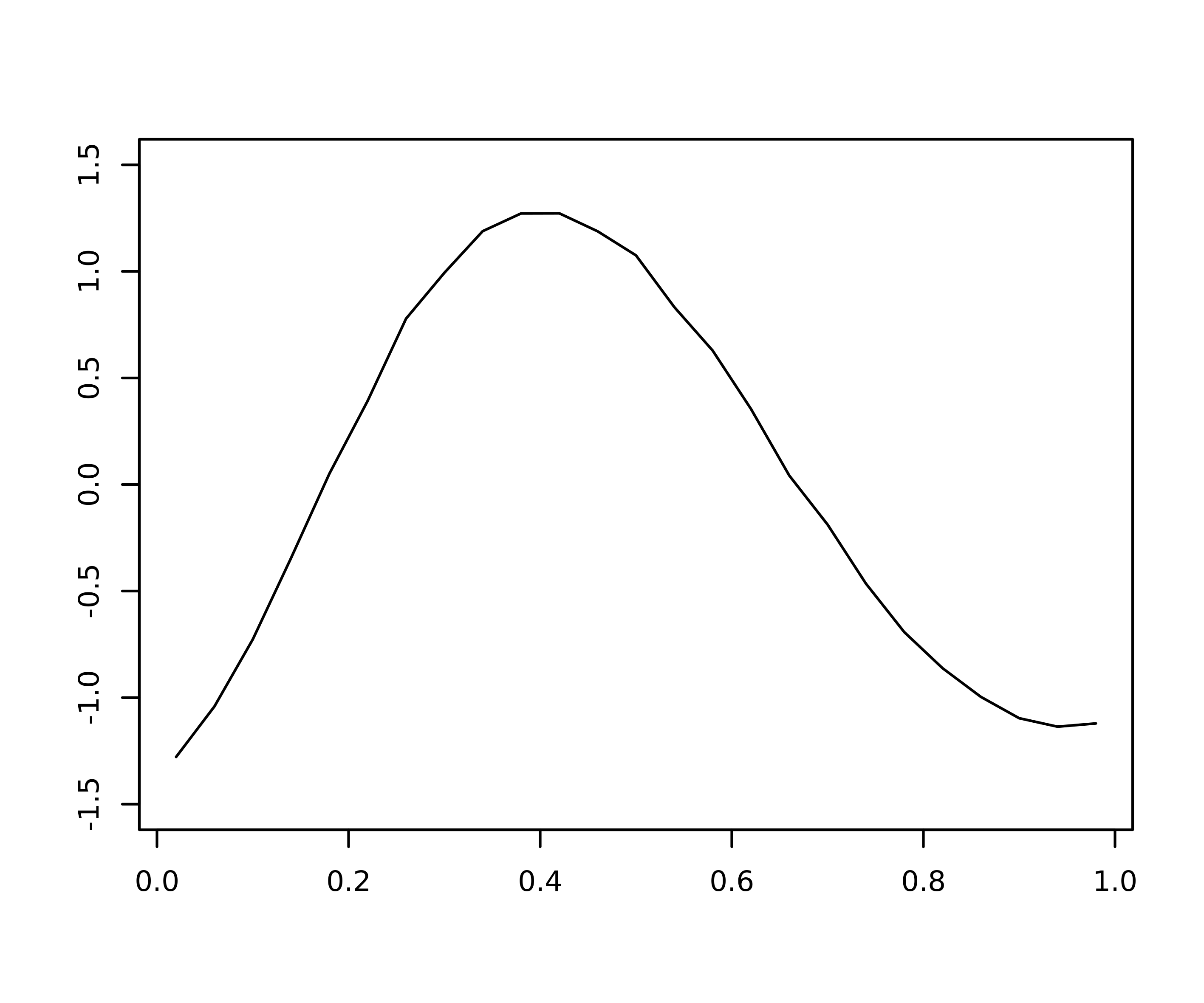}~~~
    \includegraphics[width = 6cm, height = 6 cm]{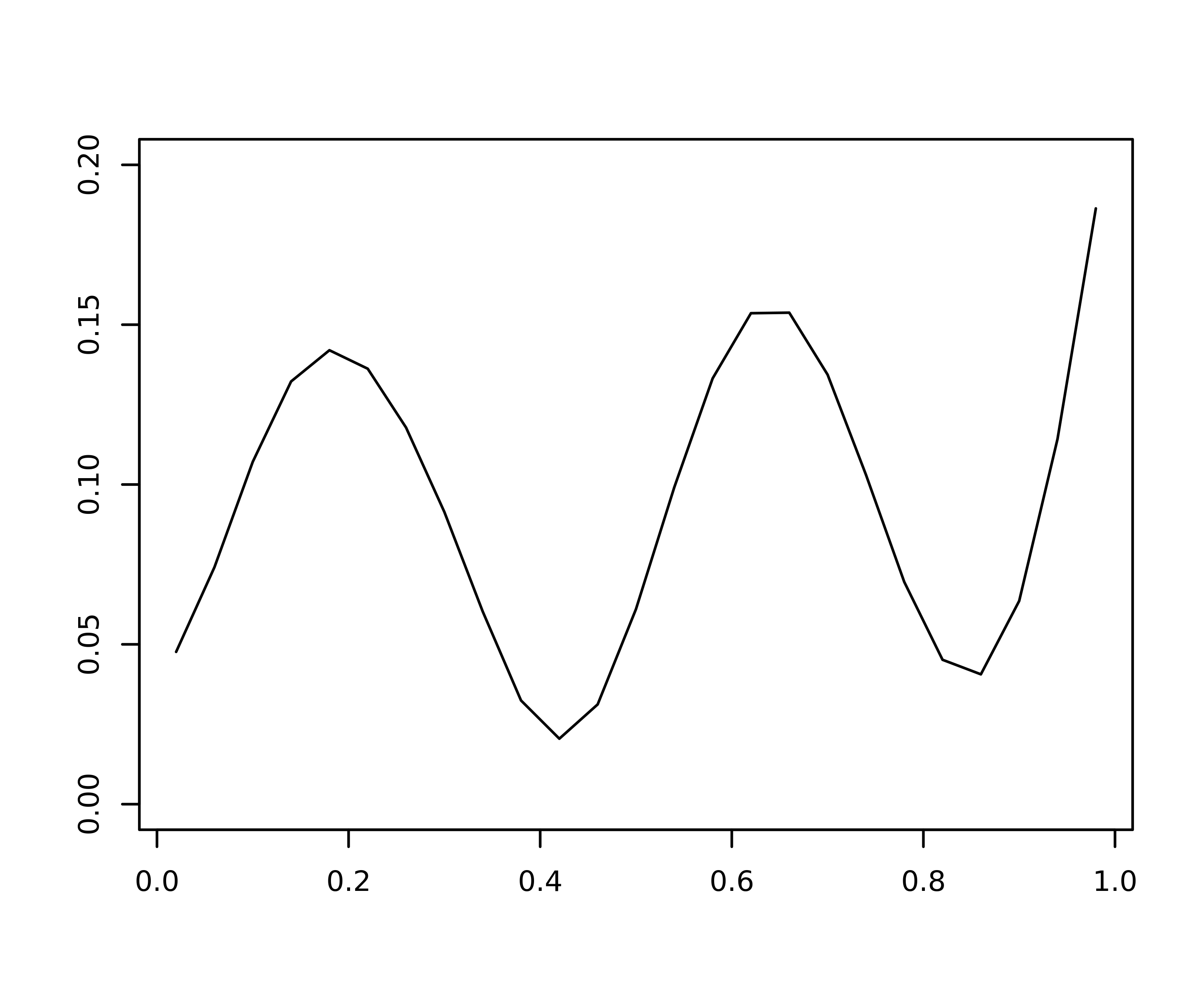}
    \caption{\it Expectation function $\mu_1$ (left panel) and variance function $\sigma_1$ (right panel)  used 
      in the simulation study in Section~\ref{sec431} and 
      Section~\ref{sec432}.   \label{fig:initial-functions}}}
\end{figure}
  
\begin{figure}[t]
  {  \centering
    \includegraphics[width = 6cm, height = 6 cm]{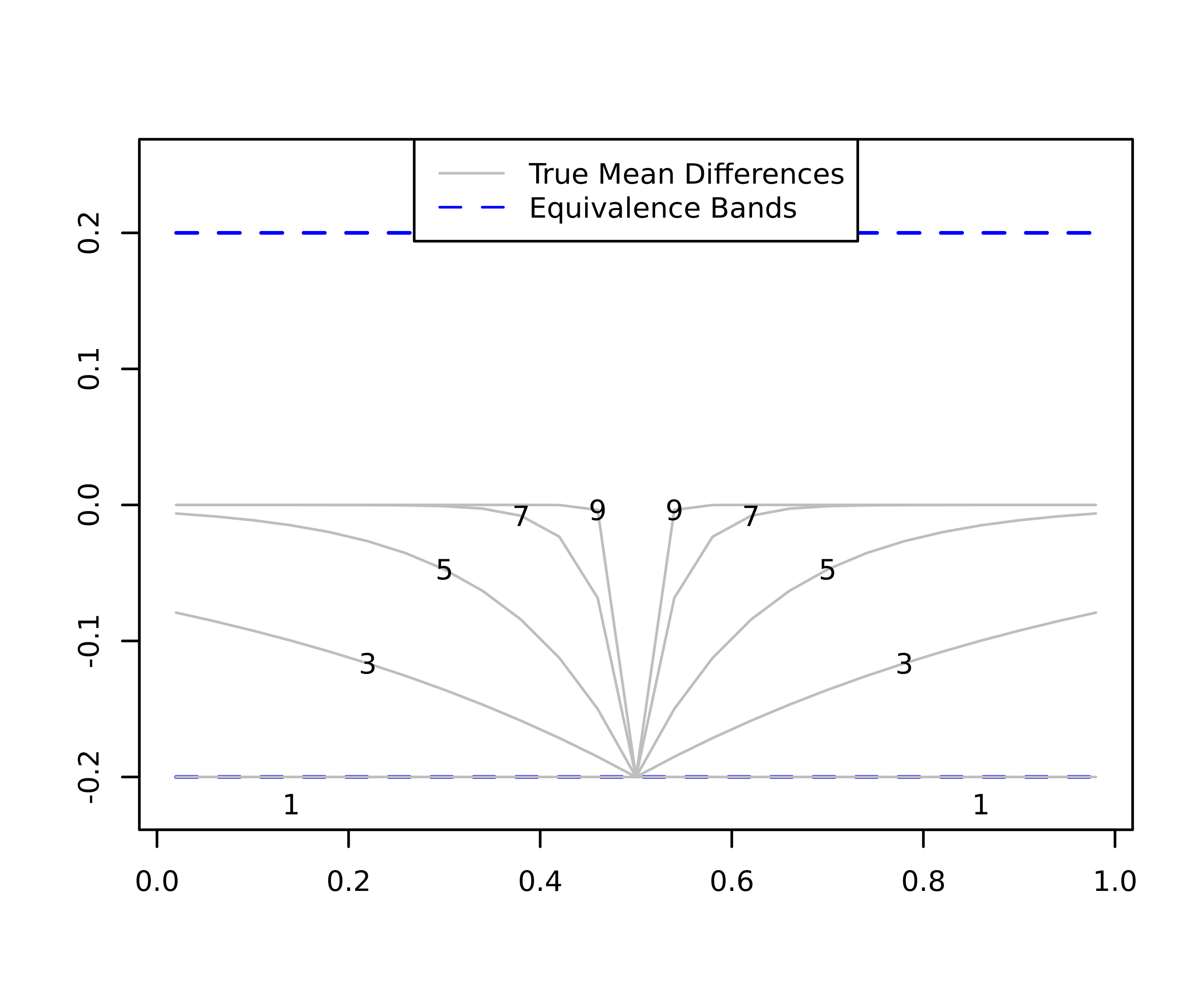}~~~
    \includegraphics[width = 6cm, height = 6 cm]{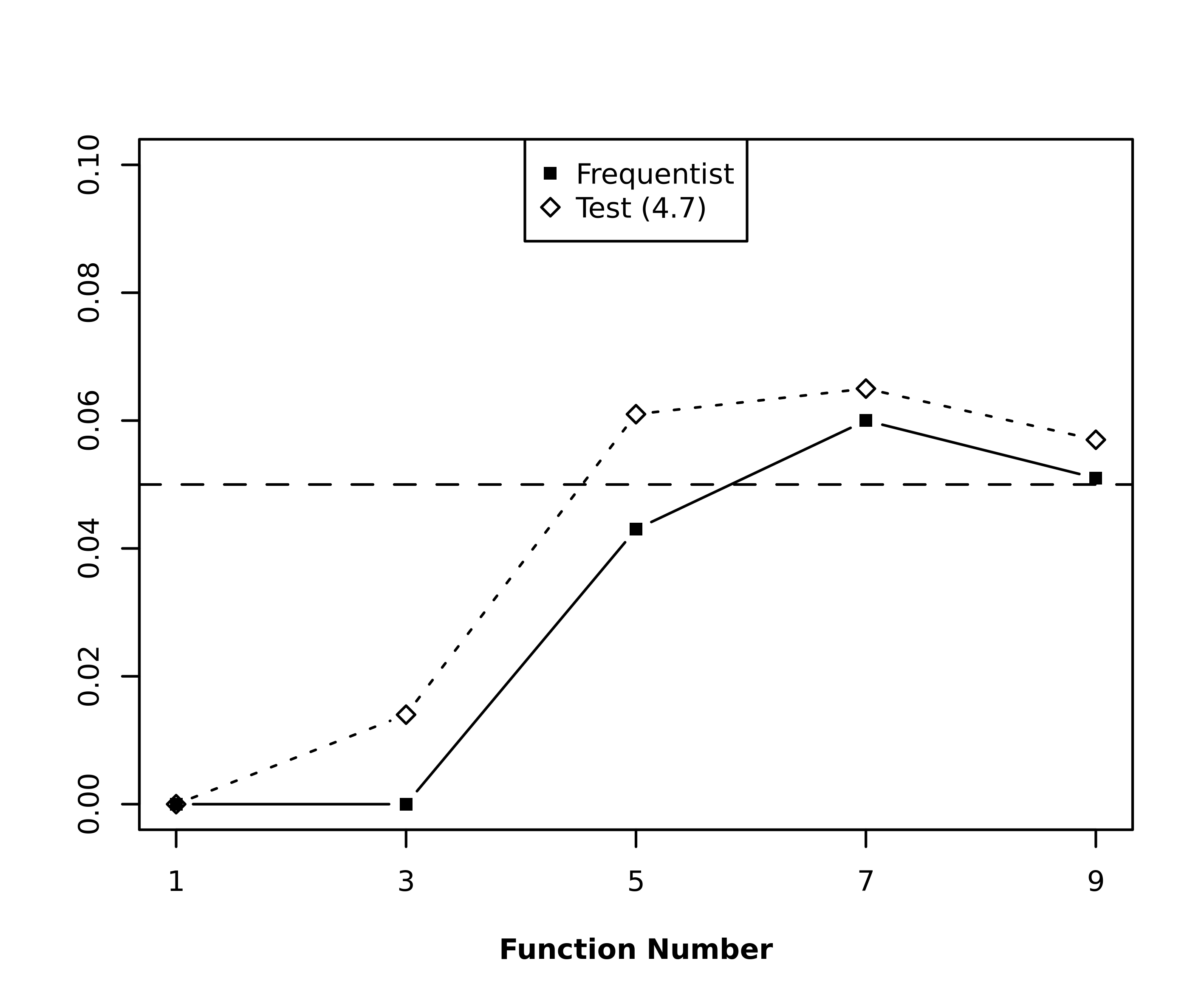}
    \caption{\it Approximation of the nominal level by  the frequentist test proposed by 
      \cite{fogarty2014} (called Frequentist) and the test \eqref{t1} for the hypotheses 
      \eqref{H0mean} with $\kappa_{l}  \equiv -0.2$ and 
      $\kappa_{u} \equiv  0.2$. 
      Left panel: True differences  for  scenarios $1$, $3$, $5$, $7$ and $9$.  
      Right part: simulated nominal level.}
    \label{fig:fogarty-size}}
\end{figure}

\begin{figure}[t]
  {  \centering
    \includegraphics[width = 6cm, height = 6 cm]{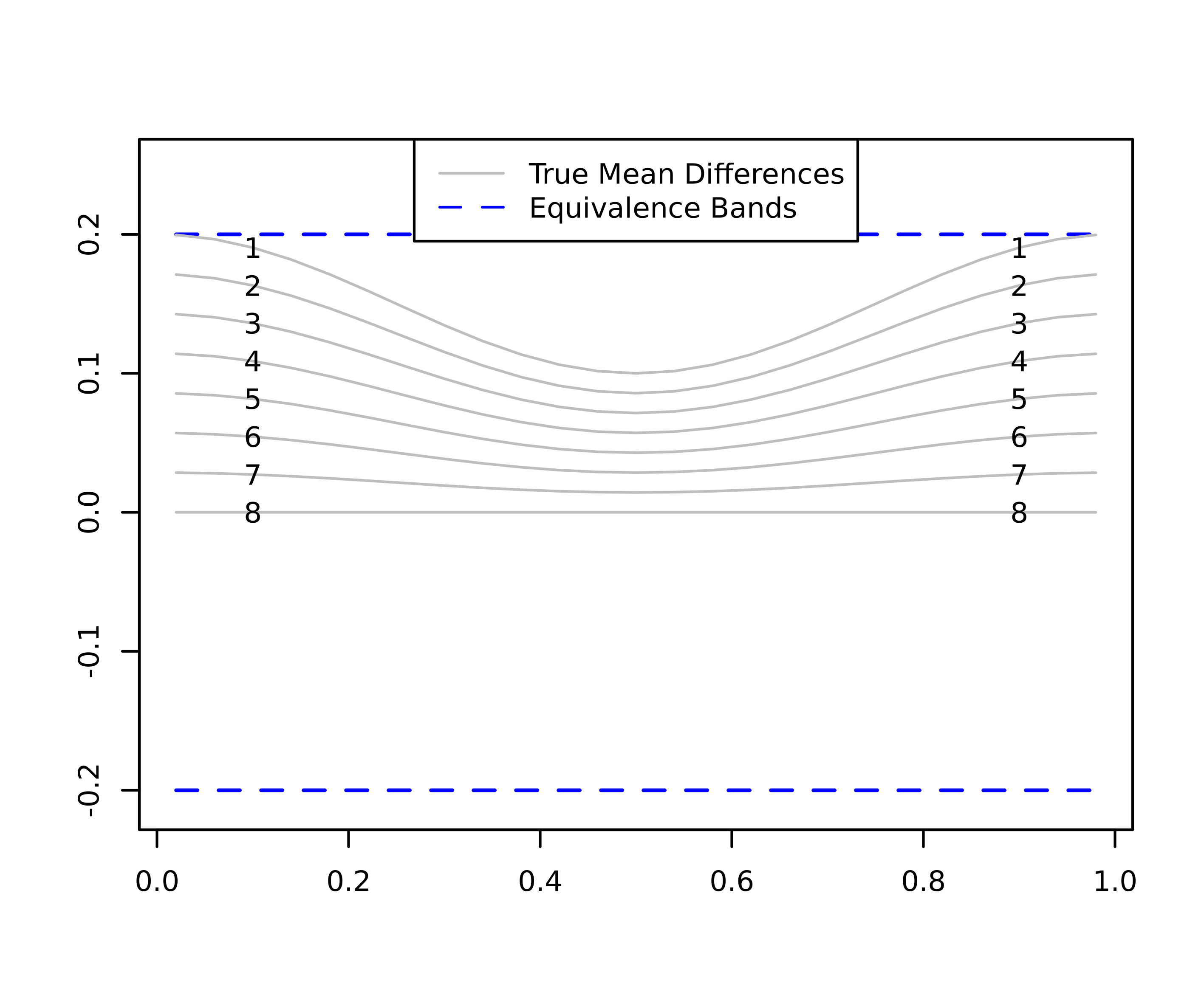}
    \includegraphics[width = 6cm, height = 6 cm]{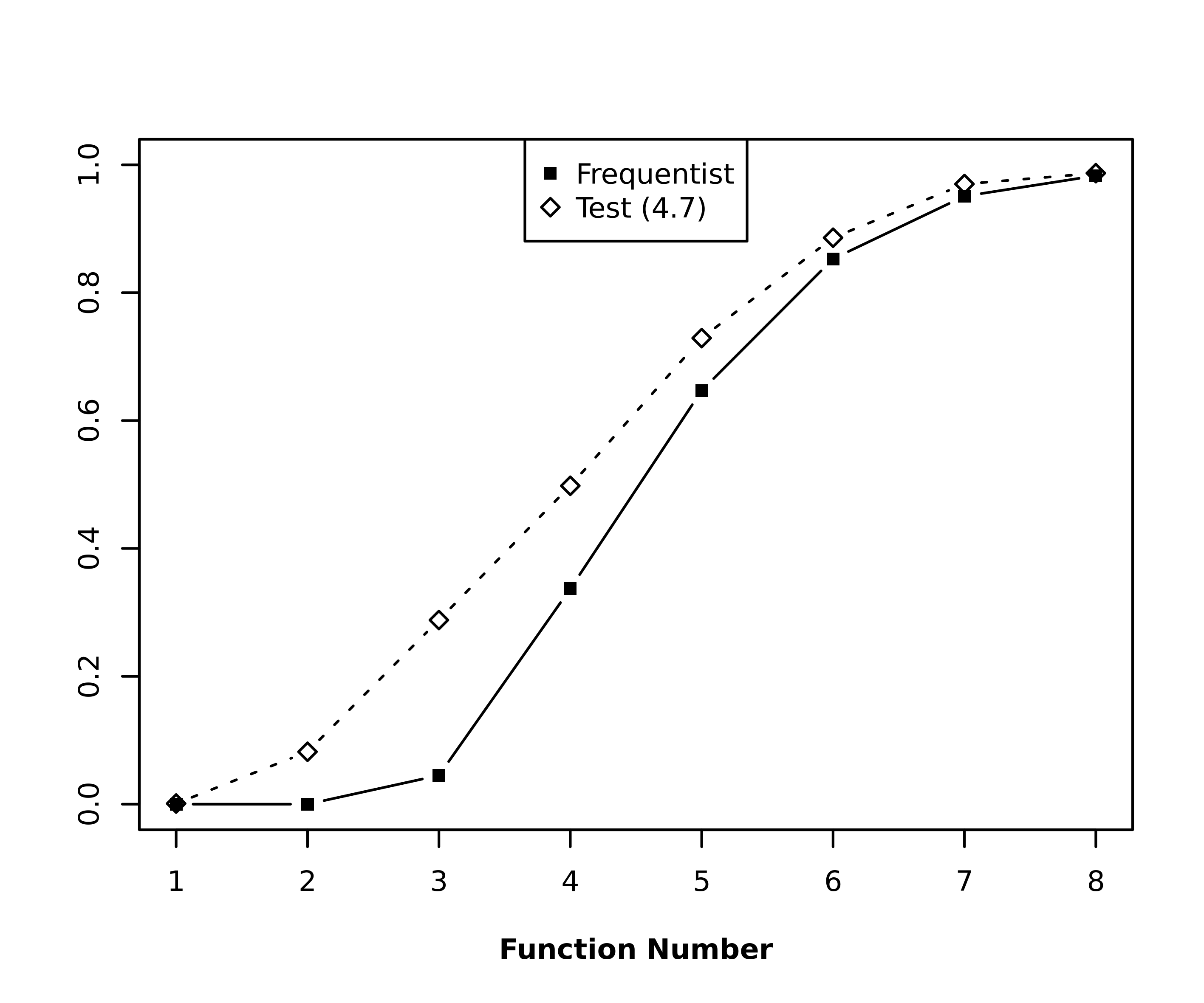}
    \caption{\it Power comparison of the frequentist test proposed by 
      \cite{fogarty2014} (called Frequentist) and the test \eqref{t1} for the hypotheses 
      \eqref{H0mean} with $\kappa_{l}  \equiv -0.2$ and 
      $\kappa_{u} \equiv  0.2$. 
Left panel: true mean difference for each scenario (1-8). Right panel: simulated rejection probabilities. 
    }
    \label{fig:fogarty-power}}
\end{figure}

\subsubsection{Comparison of mean functions}  \label{sec431}

For the mean functions, we consider five different scenarios. The mean function 
$\mu_1$ is the same in each scenario and can be obtained from the software code provided by    \cite{fogarty2014}.
It is not defined explicitly and displayed  in the left panel of 
Figure~\ref{fig:initial-functions}. 
The mean function $\mu_2$ is defined by
\begin{align*}
 \mu_{2}(t) =    \mu_{1}(t) + 0.2 \exp \big(-a_i \, |t-1/2|\big) 
\end{align*}
(thus the difference has a parametric form), 
where $a_1 = 0$ and $a_i = 10^{2(i-2)/7}$  for $i = 3,5,7,9$. 
The differences $\mu_{1 }- \mu_{2} $  correspond to the functions $1$, $3$, $5$, $7$ and $9$ in  the left part of Figure~\ref{fig:fogarty-size}, which also shows the equivalence bounds given by 
$\kappa_{l} (t) \equiv  -0.2$ and $\kappa_{u} (t) \equiv  0.2$. 
Note that 
\cite{fogarty2014} only investigate the equivalence between   the curves on the set $\{  t_j = (j-0.5)/25 \colon j = 1,\dots,25 \}$ in their simulations and for the sake of comparison, we 
consider the same set here. The  variance function 
$\sigma_1^2$ is also the same in each scenario and it is displayed in the right 
panel of Figure~\ref{fig:initial-functions}. The variance function $\sigma_2^2$ is defined by 
\begin{align}
  \frac{\sigma^2_1}{\sigma^2_2}(t) 
  = \exp\big(\log(2) \exp\big(-a_i \, |t-1/2|\big)\big) \,  \label{var1}
\end{align}
where   $i =1,3,5,7,9$.

The right part of the Figure~\ref{fig:fogarty-size} shows the simulated nominal 
level of the bootstrap test \eqref{t1} and the frequentist test proposed by  
\cite{fogarty2014} for the five cases under consideration. We observe that the 
frequentist test of \cite{fogarty2014} approximates the nominal level rather well for the function $9$, slightly exceeds the nominal level for function 7 and is conservative in the cases $1$, $3$ and $5$. The test \eqref{t1} shows a similar 
picture, where it provides a better approximation of the nominal level for the function $3$ and slightly exceeds the nominal in the cases  $5$, $7$ and  $9$.

Next we study the power of the two tests for the hypotheses \eqref{H0mean}. The 
mean function $\mu_1$ is given in the left panel of Figure~\ref{fig:initial-functions} 
and $\mu_2$ is defined by 
$$
 \mu_{2}(t) = \mu_{1 }(t) - b_i \cos(2\pi t) - c_i
$$
 $b_i = 0.05-0.1\cdot(i-1)/14$ and $c_i = 0.15-0.3\cdot(i-1)/14$ for $i=1,\dots, 8$. The variance  function
 $\sigma_1^2$ is given in  the right panel of Figure \ref{fig:initial-functions}  and $\sigma^2_2$ is defined by
\begin{align} \label{var2}
\frac{\sigma^2_1}{\sigma^2_2}(t) 
=  (0.1 \cos(2\pi t) + 1.8)^{d_i} \, ,
\end{align}
where $d_i = -1+2\cdot(i-1)/14$ for $i=1,\dots,8$.
The mean differences are depicted in the left part of Figure~\ref{fig:fogarty-power}. We observe that the frequentist test of \cite{fogarty2014} is outperformed by the new test \eqref{t1} proposed in this paper.
While the differences between the test \eqref{t1} and  the frequentist test of  \cite{fogarty2014}  are small in scenarios $6-8$ (because the power of both tests is close to $1$), 
we observe substantial advantages of the new test \eqref{t1} for the functions $2-5$.

\begin{figure}[t]
  {  \centering
    \includegraphics[width = 6cm, height = 6 cm]{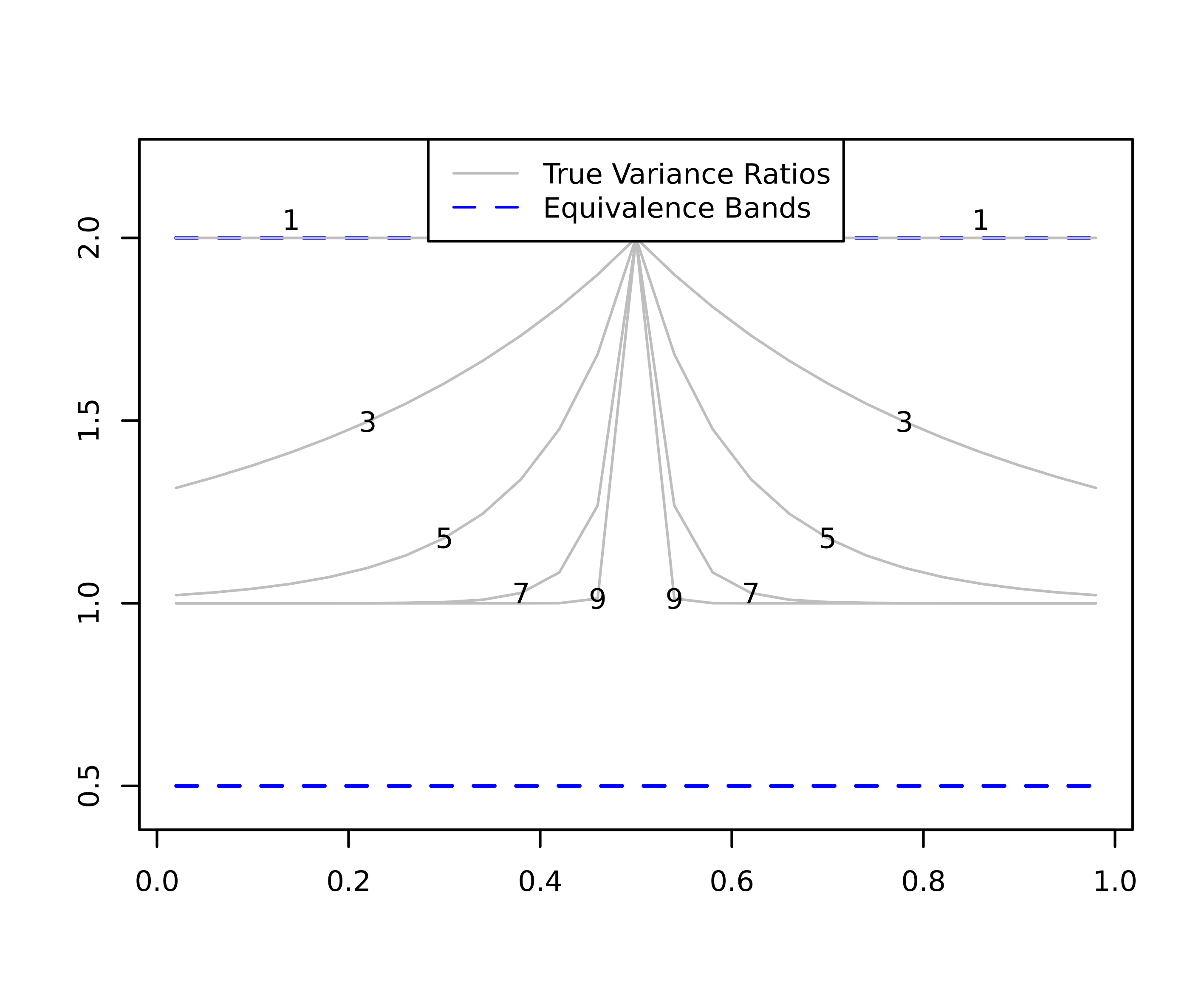}~~
    \includegraphics[width = 6cm, height = 6 cm]{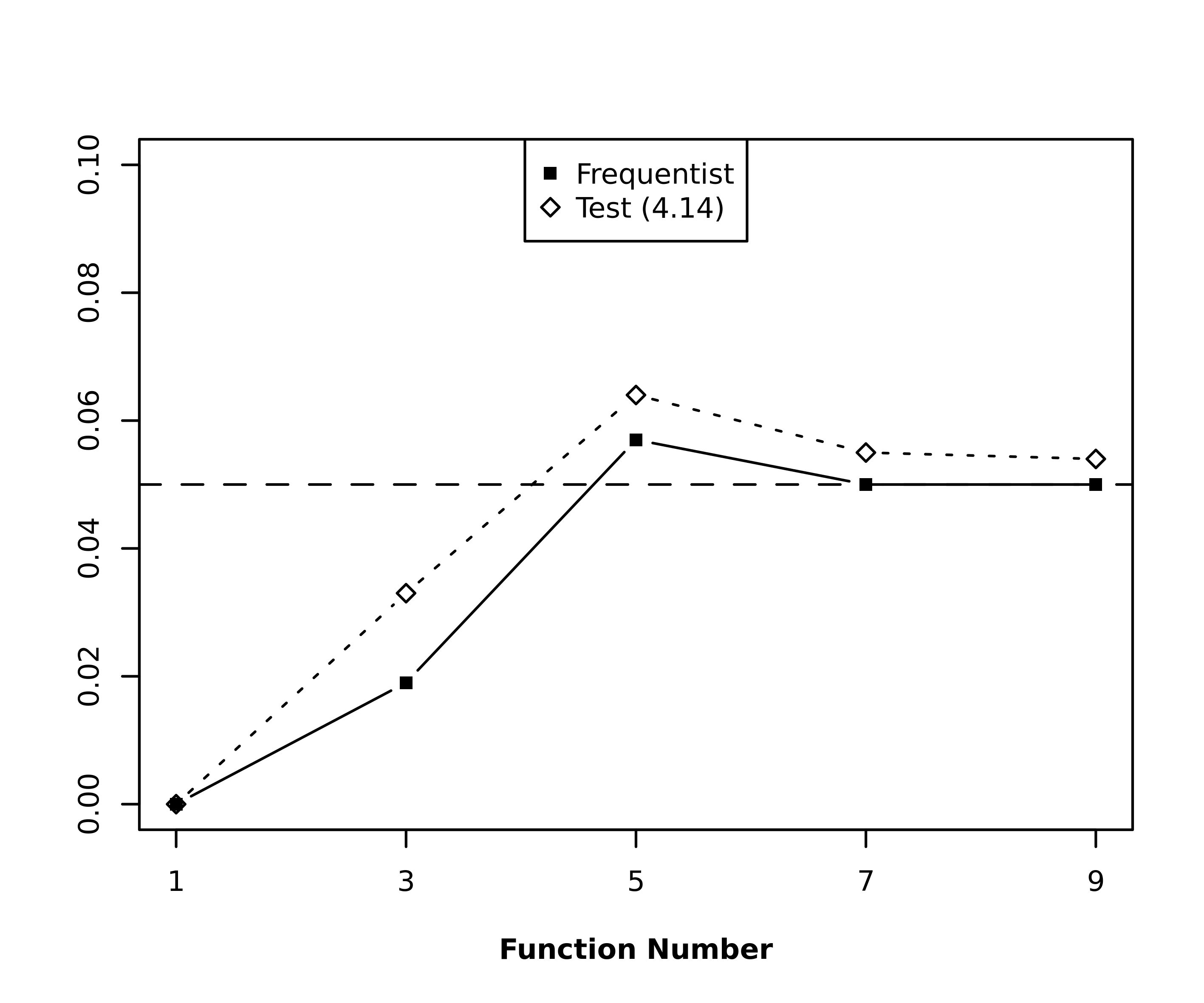}
    \caption{\it Approximation of the nominal level by  the frequentist test proposed by 
      \cite{fogarty2014} (called Frequentist) and the test \eqref{eq:var-test}
 for the hypotheses 
      \eqref{H0var} with $\zeta_l \equiv 0.5$ and $\zeta_u \equiv 2$.      Left part: True ratio of the variance functions in the scenarios 
      $1,3,5,7$ and $9$ in \eqref{var1}. 
      Right part: Simulated rejection probabilities.    }
    \label{fig:fogarty-size-var}}
\end{figure}
\begin{figure}[t]
  {  \centering
    \includegraphics[width = 6cm, height = 6 cm]{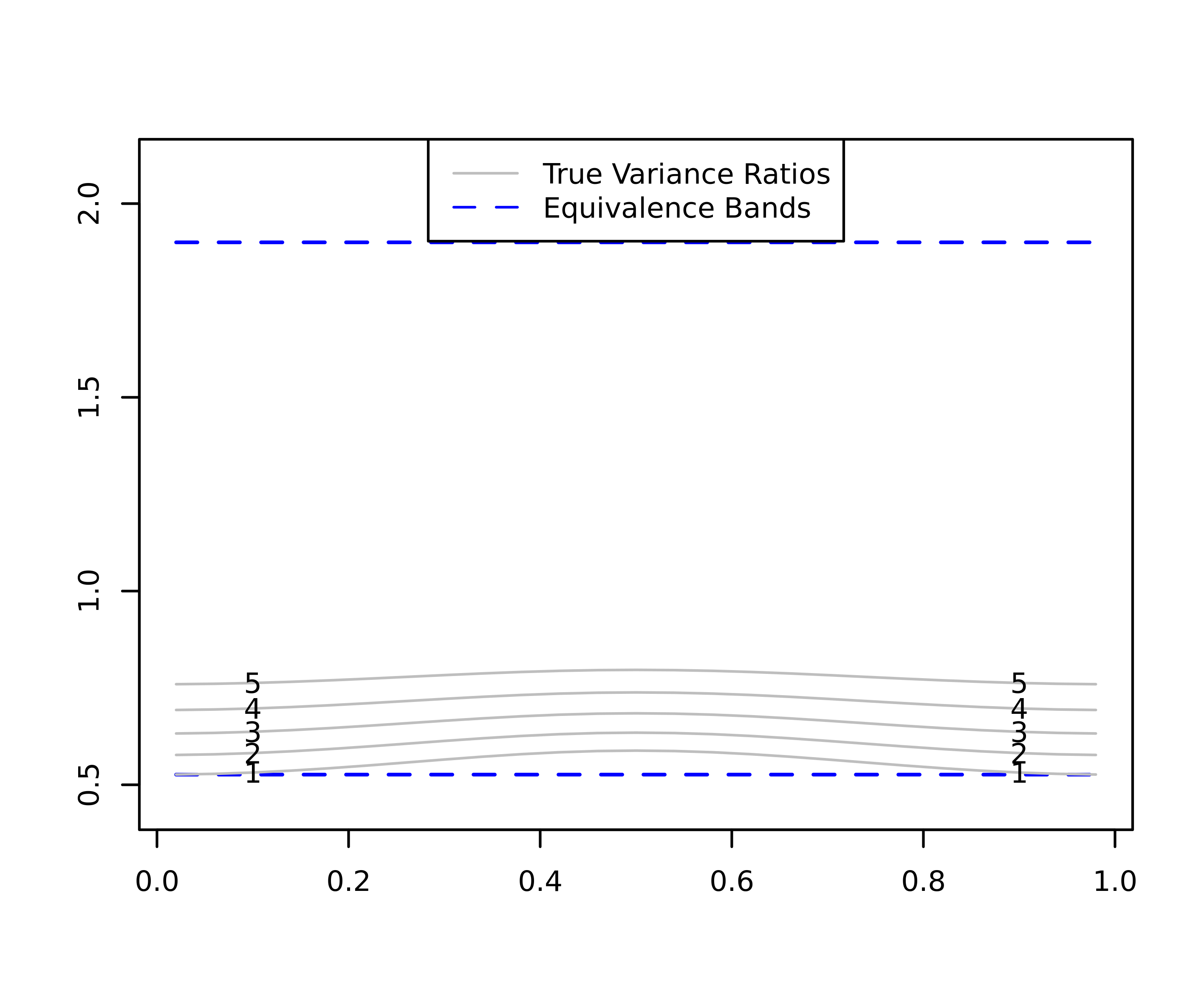}~~
    \includegraphics[width = 6cm, height = 6 cm]{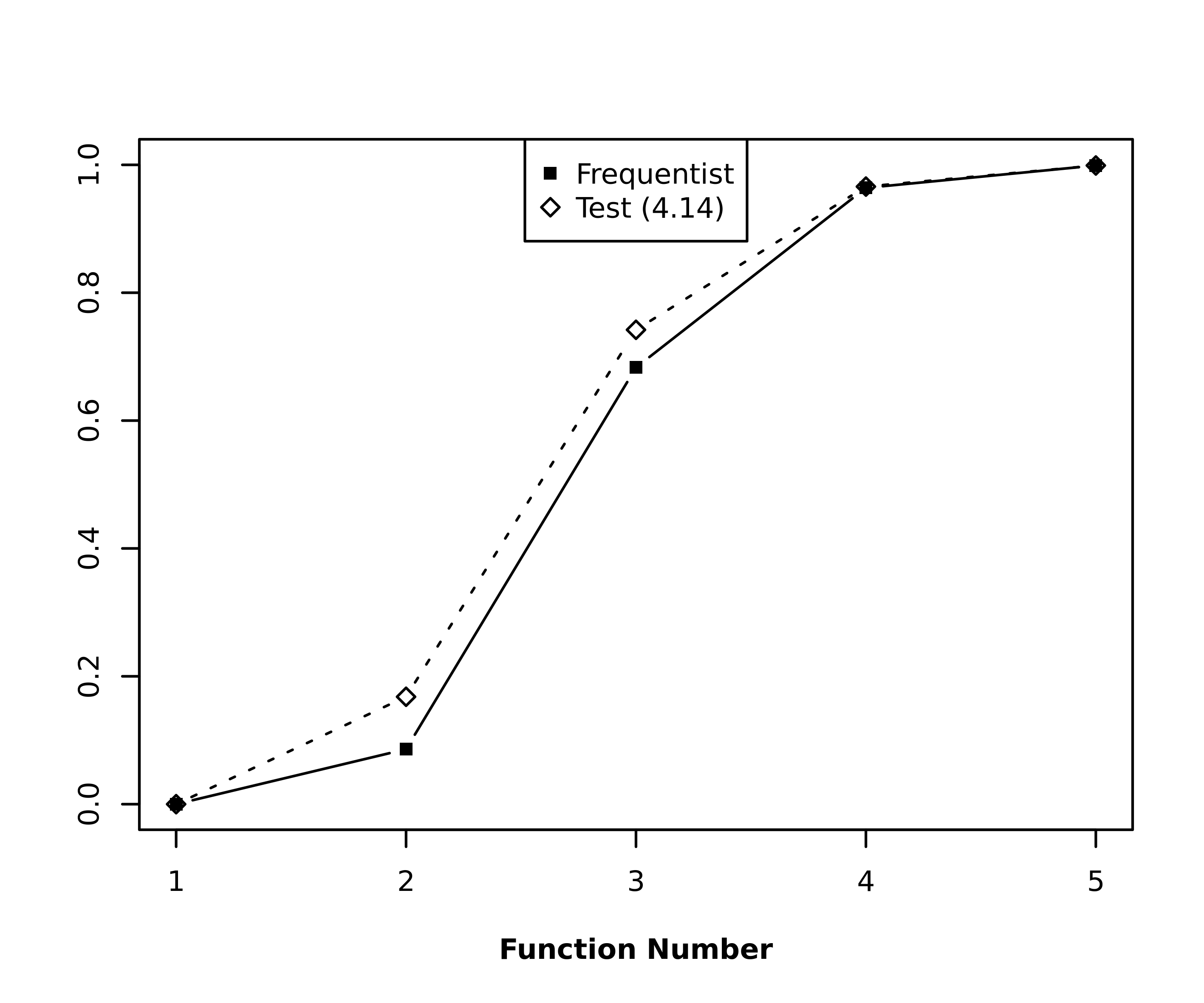}
    \caption{\it
    Power comparison of the frequentist test proposed by 
      \cite{fogarty2014} (called Frequentist) and the test \eqref{eq:var-test}  for the hypotheses 
      \eqref{H0var} with $\zeta_l \equiv 1/1.9$ and $\zeta_u \equiv 1.9$. 
      Left part: True ratio of the variance functions in the scenarios $1-5$ in 
      \eqref{var2}. 
      Right part: Empirical rejection probabilities.
    }
    \label{fig:fogarty-power-var}}
\end{figure}

\subsubsection{Variance functions} \label{sec432}

In this section, we consider the same scenarios  as in the previous section and 
investigate the finite sample properties of the tests for 
the equivalence of the variance functions of the two samples. For the different 
scenarios, the decision rule in \eqref{eq:var-test} is applied in order to 
decide for the null or the alternative hypothesis which are defined by 
\eqref{eq:equi-hypotheses-var}  or equivalently by \eqref{H0var}.    The results are then compared with those of the 
frequentist  test developed in \cite{fogarty2014}. In the left part of 
Figure~\ref{fig:fogarty-size-var}, we display the true ratio of the variance 
functions for each considered scenario in \eqref{var1}  as well as the equivalence bands defined by $\zeta_l \equiv 0.5, \zeta_u \equiv 2$. 
The right part of this figure displays the simulated nominal level of the 
bootstrap test \eqref{eq:var-test} and the frequentist test proposed by 
\cite{fogarty2014} on the boundary  of the null hypothesis. 
Similar to the results for testing the equivalence of the 
means, the frequentist approximates the test level slightly better than the 
bootstrap test in the scenarios $5, 7$ and $9$. In scenario $1$, both tests are conservative. The same is true for scenario $3$ but in this case, the empirical rejection probability of the new test is closer to the nominal level.

The true ratio of the variance functions for the considered scenarios under the 
alternative hypothesis and the used equivalence bands $\zeta_l \equiv 1/1.9, \zeta_u \equiv 1.9$
 are displayed in the left part of 
Figure~\ref{fig:fogarty-power-var}. Only the functions $1$ - $ 5$ in 
\eqref{var2} are considered since both tests always reject the null hypothesis 
in the cases $6-8$. 
The rejection rates of the two tests corresponding to the five considered 
scenarios are displayed in the right panel of Figure~\ref{fig:fogarty-power-var}. 
We observe  a superior performance of the new  bootstrap test 
\eqref{eq:var-test} in all the considered scenarios, where in the scenarios 
$1$, $4$ and $5$ the differences are very small.

\bigskip
\medskip

\noindent 
{\bf Acknowledgements}
This research was partially supported by the Collaborative Research Center `Statistical modeling
of nonlinear dynamic processes' ({\it Sonderforschungsbereich 823, Teilprojekt A1, C1})
and the Research Training Group `High-dimensional phenomena in probability - fluctuations and 
discontinuity' ({\it RTG 2131}).
The authors are grateful to Martina
 Stein, who typed  parts  of this manuscript with considerable technical expertise and to Dr. Colin Fogarty for sending us the code of the procedures developed by \cite{fogarty2014}.
 
{\small

\begin{thebibliography}{}

\bibitem[Aue et~al., 2015]{aueDubartNorinhoHormann2015}
Aue, A., {Dubart Norinho}, D., and H\"ormann, S. (2015).
\newblock On the prediction of stationary functional time series.
\newblock {\em Journal of the American Statistical Association}, 110:378--392.

\bibitem[Berger, 1982]{berger1982}
Berger, R.~L. (1982).
\newblock Multiparameter hypothesis testing and acceptance sampling.
\newblock {\em Technometrics}, 24:295--300.

\bibitem[Billingsley, 1968]{billingsley1968}
Billingsley, P. (1968).
\newblock {\em Convergence of Probability Measures}.
\newblock Wiley, New York.

\bibitem[Bradley, 2005]{bradley2005}
Bradley, R.~C. (2005).
\newblock Basic properties of strong mixing conditions. {A} survey and some
  open questions.
\newblock {\em Probability Surveys}, 2:107--144.

\bibitem[B\"ucher and Kojadinovic, 2019]{buecher2019}
B\"ucher, A. and Kojadinovic, I. (2019).
\newblock A note on conditional versus joint unconditional weak convergence in
  bootstrap consistency results.
\newblock {\em Journal of Theoretical Probability}, 32:1145--1165.

\bibitem[C\'arcamo et~al., 2020]{carcamo2019}
C\'arcamo, J., Rodr\'iguez, L.-A., and Cuevas, A. (2020).
\newblock Directional differentiability for supremum-type functionals:
  statistical applications.
\newblock {\em Bernoulli, to appear; ArXiv e-print 1902.01136}.

\bibitem[Dette et~al., 2020]{dette2018}
Dette, H., Kokot, K., and Aue, A. (2020).
\newblock Functional data analysis in the banach space of continuous functions.
\newblock {\em Annals of Statistics, to appear; ArXiv e-print 1710.07781v2}.

\bibitem[Dette et~al., 2018]{detmolvolbre2015}
Dette, H., M{\"o}llenhoff, K., Volgushev, S., and Bretz, F. (2018).
\newblock Equivalence of regression curves.
\newblock {\em Journal of the American Statistical Association}, 113:711--729.

\bibitem[Ferraty and Vieu, 2010]{FerratyVieu2010}
Ferraty, F. and Vieu, P. (2010).
\newblock {\em Nonparametric {F}unctional {D}ata {A}nalysis}.
\newblock Springer-Verlag, New York.

\bibitem[Fogarty and Small, 2014]{fogarty2014}
Fogarty, C.~B. and Small, D.~S. (2014).
\newblock Equivalence testing for functional data with an application to
  comparing pulmonary function devices.
\newblock {\em Ann. Appl. Stat.}, 8(4):2002--2026.

\bibitem[Gaenssler et~al., 2007]{Gaenssler2007}
Gaenssler, P., Moln{\'a}r, P., and Rost, D. (2007).
\newblock On continuity and strict increase of the cdf for the sup-functional
  of a gaussian process with applications to statistics.
\newblock {\em Results in Mathematics}, 51(1):51--60.

\bibitem[Gsteiger et~al., 2011]{gsteiger2011}
Gsteiger, S., Bretz, F., and Liu, W. (2011).
\newblock Simultaneous confidence bands for nonlinear regression models with
  application to population pharmacokinetic analyses.
\newblock {\em Journal of Biopharmaceutical Statistics}, 21(4):708--725.

\bibitem[Horv\'ath and Kokoszka, 2012]{HorvathKokoskza2012}
Horv\'ath, L. and Kokoszka, P. (2012).
\newblock {\em Inference for {F}unctional {D}ata with {A}pplications}.
\newblock Springer-Verlag, New York.

\bibitem[Hsing and Eubank, 2015]{hsingeubank2015}
Hsing, T. and Eubank, R. (2015).
\newblock {\em Theoretical Foundations of Functional Data Analysis, with an
  Introduction to linear Operators}.
\newblock Wiley, New York.

\bibitem[Janson and Kaijser, 2015]{jankai2015}
Janson, S. and Kaijser, S. (2015).
\newblock Higher moments of {B}anach space valued random variables.
\newblock {\em Memoirs of the American Mathematical Society}, 238.

\bibitem[Liebl and Reimherr, 2019]{lieblreim2019}
Liebl, D. and Reimherr, M. (2019).
\newblock Fast and fair simultaneous confidence bands for fundtional
  parameters.
\newblock {\em arXiv:1910.00131}.

\bibitem[Liu et~al., 2009]{liubrehaywynn2009}
Liu, W., Bretz, F., Hayter, A.~J., and Wynn, H.~P. (2009).
\newblock Assessing non-superiority, non-inferiority of equivalence when
  comparing two regression models over a restricted covariate region.
\newblock {\em Biometrics}, 65(4):1279--1287.

\bibitem[M\"ollenhoff et~al., 2018]{moedetkotvolcol2019}
M\"ollenhoff, K., Dette, H., Kotzagiorgis, E., Volgushev, S., and Collignon, O.
  (2018).
\newblock Regulatory assessment of drug dissolution profiles comparability via
  maximum deviation.
\newblock {\em Statistics in Medicine}, 37(20):2968--2981.

\bibitem[Paix{\~{a}}o et~al., 2017]{paixao2017}
Paix{\~{a}}o, P., Gouveia, L.~F., Silva, N., and Morais, J.~A. (2017).
\newblock Evaluation of dissolution profile similarity - {C}omparison between
  the {$f_2$}, the multivariate statistical distance and the {$f_2$}
  bootstrapping methods.
\newblock {\em European Journal of Pharmaceutics and Biopharmaceutics},
  79:29--50.

\bibitem[Phillips, 1990]{phillips1990}
Phillips, K.~F. (1990).
\newblock Power of the two one-sided tests procedure in bioequivalence.
\newblock {\em Journal of pharmacokinetics and biopharmaceutics},
  18(2):137--144.

\bibitem[Ramsay and Silverman, 2005]{RamsaySilverman2005}
Ramsay, J.~O. and Silverman, B.~W. (2005).
\newblock {\em Functional Data Analysis}.
\newblock Springer, New York, second edition.

\bibitem[Schuirmann, 1987]{schuirmann1987}
Schuirmann, D.~J. (1987).
\newblock A comparison of the two one-sided tests procedure and the power
  approach for assessing the equivalence of average bioavailability.
\newblock {\em Journal of pharmacokinetics and biopharmaceutics},
  15(6):657--680.

\bibitem[Van~der Vaart, 1998]{vaart_1998}
Van~der Vaart, A.~W. (1998).
\newblock {\em Asymptotic Statistics}.
\newblock Cambridge Series in Statistical and Probabilistic Mathematics.
  Cambridge University Press.

\bibitem[Van~der Vaart and Wellner, 1996]{wellner1996}
Van~der Vaart, A.~W. and Wellner, J.~A. (1996).
\newblock {\em Weak Convergence and Empirical Processes: With Applications in
  Statistics}.
\newblock Springer, New York.

\bibitem[Wellek, 2010]{wellek2010testing}
Wellek, S. (2010).
\newblock {\em Testing statistical hypotheses of equivalence and
  noninferiority}.
\newblock CRC Press.

\bibitem[Yoshida et~al., 2017]{yoshida2017}
Yoshida, H., Shibata, H., Izutsu, K.~I., and Goda, Y. (2017).
\newblock Comparison of dissolution similarity assessment methods for products
  with large variations: {$f_2$} statistics and model-independent multivariate
  confidence region procedure for dissolution profiles of multiple oral
  products.
\newblock {\em Biological and Pharmaceutical Bulletin}, 40(5):722--725.

\end{thebibliography}

}
 \newpage 

\appendix
\section{Appendix: theoretical justification} \label{sec6}
\def\theequation{5.\arabic{equation}}
\setcounter{equation}{0}

In this section we provide proofs and the necessary assumptions for our main theoretical results. We begin with some basic facts about Banach space valued random variables.

\subsection{ $C([0,1])$-valued random variables } \label{sec61}

Throughout this paper we assume that all random variables are elements of the 
space $C([0,1])$ of all continuous functions from the compact set $[0,1]$ into 
$\R$. The space $C([0,1])$ is equipped with the sup-norm defined by 
$\|f\|_\infty=\sup_{t\in [0,1]}|f(t)|$.  It is assumed that the underlying probability 
space $(\Omega,\mathcal{A},\mathbb{P})$ is complete and speaking of 
measurability is always meant with respect to the natural Borel $\sigma$-field 
$\mathcal{B}([0,1])$ (generated by the open sets relative to the sup-norm 
$\|\cdot\|_\infty$). Theorem 11.7 in \cite{jankai2015} implies that $C([0,1])$ is 
separable and measurability issues are avoided. Completeness and separability 
of $C([0,1])$ directly imply that any random variable $X$ in $C([0,1])$ is 
tight \citep[see Theorem 1.3 in][]{billingsley1968}).

Expectations and higher-order moments of $C([0,1])$-valued random variables can 
be defined formally in different ways, for example through injective tensor 
products \citep[see][]{jankai2015}. Denote by $\mathbb{E}[X]$ the expectation of 
a random variable $X$ in $C([0,1])$ and note that it exists as an element of 
$C([0,1])$ whenever $\mathbb{E}[\|X\|_\infty] < \infty$. Generally, the $k$th moment 
of $X$ exists as an element of $C([0,1]^k)$ whenever 
$\mathbb{E}[\|X\|_\infty^k] = \mathbb{E}[ \sup_{t \in [0,1]} |X(t)|^k] < \infty $ 
and it can be computed through pointwise evaluation as 
$\mathbb{E} [X (t_1) \cdots X(t_k)]$  \citep[see Chapter 11 of][]{jankai2015}). In 
particular, covariance kernels of random variables in $C([0,1])$ can be 
computed in a pointwise fashion and the variance function of $X$
can be defined by $\sigma^2(t) = \E[(X(t)-\E[X(t)])^2]$.
A  random variable $X\in C([0,1])$ is said to be Gaussian if all 
finite dimensional vectors 
$(X(t_1), \dots,X(t_k)) $  are multivariate normal distributed (for any 
$t_1,\dots,t_k\in [0,1]$ and $k\in\N$). The distribution 
of Gaussian random variables in $C([0,1])$ is completely characterized by its 
expectation and its covariance function 
\citep[see Chapter~2 of ][]{billingsley1968}.
 Throughout 
this paper, weak convergence in $C([0,1])$ is denoted by the symbol  
``$\rightsquigarrow$''  and the symbol $ {\stackrel{\mathcal{D}}{\longrightarrow}}$ denotes weak convergence of a finite dimensional 
random variable.

\subsection{Proofs  of the results in Section \ref{sec3} }  \label{sec62}

In this section we provide rigorous arguments  for the statements made in Section \ref{sec3}. In the following discussion 
$BL_1(C([0,1]))$ denotes the space of bounded (by $1$) Lipschitz functions from $C([0,1])$ 
into $\R$. That is the set of all  functions $h: C([0,1])\to\R$ with $\|h\|_\infty \leq 1$ and 
$|h(x)-h(y)| \leq \|x-y\|_\infty (= \sup_{t\in [0,1]}|x(t)-y(t)|)$ for any $x,y\in C([0,1])$.

\subsubsection{Basic assumptions and a limit theorem for the maximum deviation estimate}

For the proofs of  the results in Section  \ref{sec3} we make  the following assumption, which 
guarantees the existence of  the central limit 
theorem for independent random variables in $C([0,1])$ (see the discussion at the beginning of the proof of Theorem~\ref{thma1}).

\begin{assumption} \label{as:ts}
Let $(X_{1j}\colon j\in\mathbb{N})$ and 
$(X_{2j}\colon j\in\mathbb{N})$ denote two independent sequences of  
$C([0,1])$-valued random variables such that each sequence has independent 
identically distributed elements and assume that the following conditions are
satsified: \vspace{-.2cm}
\begin{enumerate}\itemsep-.4ex
  \item[(A1)]
  There exist   constants $\nu_1, \nu_2>0$, $K$ such that, for all $j\in\mathbb{N}$ and $i=1,2$, 
  \begin{align*}
  \mathbb{E}\big[\|X_{ij}\|_\infty^{2+\nu_i}\big] \leq K \, , 
  \quad  \mathbb{E}\big[\|X_{ij}\|_\infty^{4}\big] < \infty \, .
  \end{align*}
  \item[(A2)]
There exists a constant  $\vartheta>1/4$ and  
  a real-valued non-negative random variable $M_i$ 
  with $\E\big[M_i^4\big] < \infty$ such that, for  $i=1,2$ and  any $j\in\mathbb{N}$, 
  the inequality
  \begin{align*}
  |X_{ij}(s)-X_{ij}(t)|\leq  M_i \, \rho(s,t) = M_i \, |s-t|^\vartheta
  \end{align*}
  holds almost surely for all $s,t\in [0,1]$.
\end{enumerate}
\end{assumption}

\begin{theorem} \label{thma1} 
If Assumption \ref{as:ts} holds, we have 
\begin{align} \label{eq:limit-dist2}
\sqrt{n+m} \big ( \hat T^\theta_{m,n} - T^\theta \big) {\stackrel{\mathcal{D}}{\longrightarrow}} \,
Z_{\mathcal{E}, \theta} = \max\Big\{ 
\sup_{t\in \mathcal{E}^l_ \theta} \big(-Z(t)\big), \,
\sup_{t\in \mathcal{E}^u_ \theta} Z(t) \Big\} \, , 
\end{align}
where   $T^\theta$ and $\hat T^\theta_{m,n}$ are defined by  \eqref{eq:equi-hypotheses} and \eqref{eq:statistic}, respectively
and the   extremal sets $\mathcal{E}^l_ \theta$ and $\mathcal{E}^u_ \theta$ are defined by \eqref{eq:sets}. 
\end{theorem}

\noindent 
\textbf{Proof.} 
Note that condition (A2) and the fact that $\vartheta>1/4$ imply 
$\int_0^{\tilde\tau}  D(\omega,\rho)^{1/4}~d\omega <\infty$
for some $\tilde\tau > 0$ where $D(\omega,\rho)$ denotes the packing number 
with respect to the metric $\rho(s,t) = |s-t|^\vartheta$ that is the maximal 
number of $\omega$-seperated points in $[0,1]$ \citep[see][]{wellner1996}. 
Therefore  it follows 
from Theorem~2.1 in \cite{dette2018} that 
\begin{align} \label{eq:CLT}
  \hat{Z}_{m,n} = \sqrt{m+n} \, (\hat{\theta}_{m,n} - \theta) \rightsquigarrow Z
\end{align}
in $C([0,1])$ as $m,n\to\infty$ where $Z$ is a (tight) centred Gaussian random 
function in $C([0,1])$ with covariance kernel as defined in \eqref{eq:Z-kernel} \citep[see also Remark~2.1 (b) in ][]{dette2018}.
Observing the estimate 
\begin{align*}
&\sqrt{m+n} \big ( \hat T^\theta_{m,n} - T^\theta \big)   \\
 &= \sqrt{m+n} \Big( \max\Big\{ 
 \sup_{t\in [0,1]} \big(-\hat{\theta}_{m,n}(t) + \kappa_l(t) \big), \,
 \sup_{t\in [0,1]} \big(\hat{\theta}_{m,n}(t) - \kappa_u(t) \big) \Big\} 
 - T^\theta \Big) \\
 &= \sqrt{m+n} \Big( \max\Big\{ 
 \sup_{t\in \mathcal{E}_\theta^l} \big(-\hat{\theta}_{m,n}(t) + \kappa_l(t) 
 \big), \,
 \sup_{t\in \mathcal{E}_\theta^u} \big(\hat{\theta}_{m,n}(t) - \kappa_u(t) 
 \big) \Big\} 
 - T^\theta \Big) + o_\P(1) \\
 &= \sqrt{m+n} \max\Big\{ 
 \sup_{t\in \mathcal{E}_\theta^l} \big(-\hat{\theta}_{m,n}(t) + \kappa_l(t) 
 - T^\theta \big), \,
 \sup_{t\in \mathcal{E}_\theta^u} \big(\hat{\theta}_{m,n}(t) - \kappa_u(t) 
 - T^\theta \big) \Big\} + o_\P(1) \\
 &= \sqrt{m+n} \max\Big\{ 
 \sup_{t\in \mathcal{E}_\theta^l} \big(-\hat{\theta}_{m,n}(t)+\theta(t) \big), \,
 \sup_{t\in \mathcal{E}_\theta^u} \big(\hat{\theta}_{m,n}(t) - \theta(t) \big) \Big\} + o_\P(1) \, ,
\end{align*} 
the assertion of Theorem \ref{thma1} follows from \eqref{eq:CLT} and the continuous mapping theorem.
\hfill $\Box$

\subsubsection{Proof of Theorem \ref{thm1}}  \label{sec622}

We begin showing that the bootstrap process $\hat{Z}_{m,n}^{*(r)}$ converges 
conditionally given the data $(X_{lj} \colon j\in \N; ~ l=1,2)$, 
to the same limit as $\hat{Z}_{m,n}$. More precisely, this means  
\begin{align} \label{eq:conditional-convergence}
  \sup_{h\in BL_1(C([0,1]))} | \E_M h(\hat{Z}_{m,n}^{*(r)}) - \E h(Z)| = o_\P(1)
\end{align}
as $m,n\to\infty$ where $\E_M$ denotes the conditional expectation given the 
 data $(X_{lj} \colon j\in \N; ~ l=1,2)$ and the random variable $Z$ is defined by \eqref{eq:Z-kernel} 
 \citep[see for example Section~23.2.1 in][]{vaart_1998}. 

Note that the convergence in \eqref{eq:conditional-convergence} holds under the 
null and under the alternative hypothesis. 
By the continuous mapping theorem and similar arguments as given in Lemma~B.3 of \cite{dette2018} it follows that the bootstrap statistic 
$\hat T^{\theta,*(r)}_{m,n}$ defined by \eqref{hda} converges conditionally 
given the  data $(X_{lj} \colon j\in \N; ~ l=1,2)$ to the same limit as 
$\sqrt{m+n} \, (\hat T^\theta_{m,n} - T^\theta)$ that is $Z_{\cal E,\theta}$ 
(see \eqref{eq:limit-dist}).
If $T^\theta = 0$, Lemma~4.2 in \cite{buecher2019} directly implies the first 
assertion of Theorem~\ref{thm1} that is
\begin{align*}
   \lim_{m,n,R\to\infty} \mathbb{P}\big( 
  \sqrt{m+n} \, \hat{T}_{m,n}^\theta < z_{m,n,\alpha}^{*(R)} \big) = \alpha 
\end{align*}
(note that the continuity of the random variable $Z_{\cal E,\theta}$ is implied 
by the results in \cite{Gaenssler2007}).
If $T^\theta \neq 0$, write 
\begin{align*}
  \mathbb{P}\big(\sqrt{m+n} \, \hat{T}_{m,n}^\theta < z_{m,n,\alpha}^{*(R)}\big)
  = \mathbb{P}\big( 
  \sqrt{m+n} \, (\hat{T}_{m,n}^\theta - T^\theta) + \sqrt{m+n} \, T^\theta 
  < z_{m,n,\alpha}^{*(R)} \big) \, .
\end{align*}
Then it follows from \eqref{eq:limit-dist2}, \eqref{eq:conditional-convergence} and simple arguments that, for any 
$R\in\N$,
\begin{align*}
\lim_{m,n\to\infty} \mathbb{P}\big( 
  \sqrt{m+n} \, \hat{T}_{m,n}^\theta < z_{m,n,\alpha}^{*(R)} \big) = 0 ~~~ \text{and} ~~~ \liminf_{m,n\to\infty} \mathbb{P}\big(\sqrt{m+n} \, \hat{T}_{m,n}^\theta < z_{m,n,\alpha}^{*(R)} \big) = 1
\end{align*}
if $T^\theta > 0$ and $T^\theta \leq 0$, respectively. This
proves the remaining assertions of Theorem~\ref{thm1}.

In order to prove the convergence in \eqref{eq:conditional-convergence}, we 
will utilize the link between weak convergence in the Banach space of 
continuous functions $C([0,1])$ and weak convergence of empirical processes in 
the space of bounded functions $l^\infty(\F)$ from an appropriately defined   function space $\F$
 into $\R$. In fact we use the CLT derived in the proof of Theorem \ref{thma1}  (see equation \eqref{eq:CLT}) and conclude that 
Theorem~23.7 in \cite{vaart_1998} can be applied to show weak convergence of 
the empirical bootstrap process in $l^\infty(\cal F)$ conditionally given the 
data $(X_{lj} \colon j\in \N; ~ l=1,2)$. Afterwards it will be argued that this again implies the 
convergence in \eqref{eq:conditional-convergence}. 

For any $t\in [0,1]$ and $x\in C([0,1])$, consider the canonical projection $\pi_t\colon C([0,1]) \to \R$ with $\pi_t(x) = x(t)$, define the function class
\begin{align*}
  \mathcal{F} = \{ \pi_{t} \colon t \in [0,1] \}
\end{align*}
and note that this class is a subset of $C([0,1])^*$, the dual space of 
$C([0,1])$. Defining the map $x^{**}:\F\to\R$ by 
$x^{**}(\pi_t) = \pi_t(x) = x(t)$, for any $x\in C([0,1])$, leads to an 
isometric identification of $C([0,1])$ with a subset of $l^\infty(\F)$ and in 
the following, this subset is denoted by $C([0,1])^{**}$. Both, $C([0,1])$ and 
$l^\infty(\F)$, are equipped with the respective sup-norm $\|\cdot\|_\infty$ 
and it is clear that $\|x\|_\infty = \|x^{**}\|_\infty$. For any 
$C([0,1])$-valued random variable $X$ the corresponding random variable 
$X^{**}\in C([0,1])^{**} \subset l^\infty(\F)$ is defined by 
$$
  X^{**}(\pi_t) = \delta_X \pi_t 
  = \int_{C([0,1])} \pi_t(x) \delta_X(dx) = X(t)
$$ 
where $\delta_X$ denotes the dirac measure. Next, show that the weak 
convergence of a sequence of random variables in $C([0,1])$, that is 
$X_n \rightsquigarrow X$, is equivalent to the weak convergence 
$X_n^{**} \rightsquigarrow X^{**}$ in $l^\infty(\F)$  
\citep[see also Section~2.1.4 in][]{wellner1996}. Following Section~1.12 in \cite{wellner1996} (note 
that each random variable $X\in C([0,1])$ is separable), weak convergence of a 
sequence $X_n\subset C([0,1])$ to a separable random variable $X\in C([0,1])$, 
denoted by $X_n \rightsquigarrow X$, is equivalent to 
\begin{align*}
\sup_{h\in BL_1(C([0,1]))} | \E^* h(X_n) - \E h(X)| = o(1) \, .
\end{align*}
Each function $h\in BL_1(C([0,1]))$ can be identified by the function 
$h^{**}:C([0,1])^{**} \to \R$ defined through $h^{**}(x^{**}) = h(x)$ for any 
$x^{**}\in C([0,1])^{**}$. Note that $\|h^{**}\|_\infty = \|h\|_\infty \leq 1$ 
and 
$$
|h^{**}(x^{**})-h^{**}(y^{**})| = |h(x)-h(y)| \leq \|x-y\|_\infty = \|x^{**}-y^{**}\|_\infty.
$$ 
Thus $h^{**}\in BL_1(C([0,1])^{**})$ and 
\begin{align} \label{eq:BL-equality}
\begin{split}
  o(1) &= \sup_{h\in BL_1(C([0,1]))} |\E^* h(X_n) - \E h(X)| \\ 
  &= \sup_{h^{**}\in BL_1(C([0,1])^{**})} 
  |\E^* h^{**}(X_n^{**})-\E h^{**}(X^{**})| \\
  &= \sup_{g\in BL_1(l^\infty(\F))} 
  |\E^* g|_{C([0,1])^{**}}(X_n^{**})-\E g|_{C([0,1])^{**}}(X^{**})| \\
  &= \sup_{g\in BL_1(l^\infty(\F))} |\E^* g(X_n^{**})-\E g(X^{**})|
\end{split}
\end{align}
where $g|_{C([0,1])^{**}}$ denotes the restriction of the function $g$ to 
$C([0,1])^{**}$. Consequently $X_n^{**} \rightsquigarrow X^{**}$ in 
$l^\infty(\F)$ (by construction, $X$ and $X^{**}$ can be defined on the same 
original probability space and we have that $X^{**}$ is separable if and only 
if $X$ is separable).

Since the envelope function of $\F$, $F(x)=\|x\|_\infty$, is finite for any 
$x\in C([0,1])$ and the convergence in \eqref{eq:CLT} together with the 
previous discussion means that $\F$ is a Donsker class, Theorem~23.7 in 
\cite{vaart_1998} can be applied. For this purpose, note that 
$$
\hat Z_{m,n}^{*(r)} = \frac{\sqrt{m+n}}{\sqrt{m}} \hat Z_{1,m}^{*(r)} + \frac{\sqrt{m+n}}{\sqrt{n}} \hat Z_{2,n}^{*(r)},
$$ where
\begin{align*}
  \hat Z_{1,m}^{*(r)} =
  \frac {1}{\sqrt{m}} \sum^m_{j=1} (X^{*(r)}_{1j} - \overline{X}_{1 \cdot}) 
  &= \frac {1}{\sqrt{m}} \sum^m_{j=1} (M_{1j}^{(r)} - 1) X_{1j}
\end{align*}
and $\hat Z_{2,n}^{*(r)}$ is defined analogously. 
The random vector 
$(M_{11}^{(r)},\dots,M_{1m}^{(r)})$ follows a multinomial distribution with 
parameters $m$, $p_j = 1/m$, $j = 1,\dots,m$. The corresponding empirical 
bootstrap process can be written as 
\begin{align*}
  (\hat Z_{1,m}^{*(r)})^{**} 
  = \frac{1}{\sqrt{m}} \sum^m_{j=1} (M_{1j}^{(r)} - 1) \delta_{X_{1j}} \, .
\end{align*}
Now, Theorem~23.6 in \cite{vaart_1998} implies that 
\begin{align*}
  \sup_{g\in BL_1(l^\infty(\F))} 
  |\E_M g((\hat Z_{1,m}^{*(r)})^{**})-\E g(Z_1^{**})| 
  = o_\P(1) 
\end{align*}
where $\E_M$ denotes the conditional expectation given the  data $(X_{lj} \colon j\in \N; ~ l=1,2)$  and 
$Z_1^{**}$ is the (unconditional) limit of 
$\sqrt{m}(\overline{X}_{1\cdot}-\mu_1)^{**}$. Similar arguments as in \eqref{eq:BL-equality} and the subsequent discussion yield that this equation  is equivalent to 
\begin{align*}
\sup_{h\in BL_1(C([0,1]))} |\E_M h(\hat Z_{1,m}^{*(r)})-\E h(Z_1)| = o_\P(1) 
\end{align*}
which means that the sequence $\hat Z_{1,m}^{*(r)}$ converges conditionally 
given the data $(X_{lj} \colon j\in \N; ~ l=1,2)$  to $Z_1$ in $C([0,1])$. The corresponding statement for 
$\hat Z_{2,n}^{*(r)}$ can be derived similarly.
Since $\hat Z_{1,m}^{*(r)}$ and $\hat Z_{2,n}^{*(r)}$ are independent  \eqref{eq:conditional-convergence} now follows, which completes the proof of  Theorem \ref{thm1}.

\subsubsection{Proof of  Remark \ref{rem3}}  \label{sec623}

In this section we consider the case of dependent data and  give some arguments why   the decision rule in \eqref{eq:equi-testdep} based on the 
multiplier bootstrap process defined by \eqref{2bProcess} yields a consistent 
and asymptotic level $\alpha$-test. For that consider the dependency 
concept of $\varphi$-mixing (see for example \cite{bradley2005}). Denote by 
$\mathbb{P}(G|F)$ the conditional probability of $G$ given $F$ and, for any two 
$\sigma$-fields $\mathcal{F}$ and $\mathcal{G}$, define
\[
  \phi(\mathcal{F},\mathcal{G})
  = \sup \big\{ |\mathbb{P}(G|F) - \mathbb{P}(G)| \colon 
  F\in \mathcal{F}, ~G\in\mathcal{G}, ~\mathbb{P}(F)>0 \big\} \, .
\]
For a given stationary sequence $(\eta_{j}\colon j\in\mathbb{N})$ of random variables in 
$C([0,1])$, denote by $\mathcal{F}^{k^\prime}_{k}$ the $\sigma$-field generated 
by $(\eta_{j}\colon k\leq j \leq k^\prime)$. Then, the $k$th 
$\varphi$\textit{-mixing coefficient} of  $(\eta_{j}\colon j\in\mathbb{N})$ is 
defined by
\[
\varphi (k)
= \sup_{k^\prime \in \mathbb{N}} \phi (\mathcal{F}_{1}^{k^\prime} ,\mathcal{F}_{k^\prime+k}^\infty)
\]
and the stationary  time series $(\eta_{j}\colon j\in\mathbb{N})$ is called $\varphi$\textit{-mixing} whenever the sequence of mixing coefficients converges to zero as $k\to\infty$.

The statement in Remark \ref{rem3} is correct, if the   following assumptions are satisfied:
\begin{enumerate}
\item[(B1)] $(X_{1,j}\colon j\in\mathbb{N})$ and $(X_{2,j}\colon j\in\mathbb{N})$ are independent stationary  time series 
 satisfying conditions  (A1) and (A2) in Assumption~\ref{as:ts}.
\item[(B2)] Both sequences are $\varphi$-mixing and the mixing coefficients satisfy, for $i=1,2$, 
\begin{align*}
\sum_{k=1}^\infty k^{1/(1/2-\bar\tau_i) } \varphi_i(k)^{1/2}  <\infty \, , ~~
\sum_{k=1}^\infty (k+1)^{}\varphi_i(k)^{1/4}  < \infty \, ,
\end{align*}
for some $\bar\tau_i \in (1/(2+2\nu_i), 1/2 )$ where the constant $\nu_i$ is the same as in (A1).
\item[(B3)] 
The window parameters $l_1,l_2$ in the definition of the bootstrap processes in \eqref{2bProcess} are defined by $l_1 = m^{\beta_1}$, $l_2 = n^{\beta_2}$ such that 
$$
0<\beta_i<\nu_i/(2+\nu_i) \, ,~~  \bar{\tau}_i > (\beta_{i}( 2 + \nu_{i}) + 1) / (2 + 2\nu_{i}) \, , 
$$
and the constants $\nu_i$ and $\bar{\tau}_i$ are given in (A1) and (B2), respectively ($i=1,2$).

\end{enumerate}

  Under these assumptions  
it follows from Theorem~2.1 in \cite{dette2018} that the CLT in \eqref{eq:CLT} also holds 
in the dependent case, where the limiting process $Z$ is a Gaussian process 
with covariance kernel 
\begin{equation*}
  C(s,t) = \frac{1}{\tau} C_1(s,t) + \frac{1}{1-\tau} C_2(s,t) 
\end{equation*}  
and $C_i(s,t) = \sum_{j=-\infty}^\infty \cov(X_{i0}(s),X_{ij}(t)) $  ($i = 1,2$).
 Similarly, by Theorem~3.3 in the same 
reference the bootstrap process \eqref{2bProcess} satisfies 
\begin{align*}
  (\hat{Z}_{m,n}, \hat{Z}_{m,n}^{**(1)},\dots,\hat{Z}_{m,n}^{**(R)})
  \rightsquigarrow (Z, Z^{(1)}, \dots Z^{(R)})
\end{align*}
in $C([0,1])^{R+1}$ as $m,n\to\infty$ where $Z^{(1)}, \dots Z^{(R)}$ are 
independent copies of $Z$. Note that this is equivalent to the corresponding 
statement in \eqref{eq:conditional-convergence} 
\citep[by Lemma~2.2 in ][]{buecher2019} and therefore,
the statement of Remark~\ref{rem3} now follows by similar arguments as given in the proof of Theorem~\ref{thm1}.

\subsection{Proofs  of the results in Section \ref{sec4} }  \label{sec624}

\subsubsection{Proofs  of the results in Section \ref{sec51} } \label{sec51-assumptions} 

Recall the model defined in \eqref{eq:fogarty-datamodel} and assume the 
following: 
\begin{assumption} \label{(C)}
 The differences of the error terms are sampled from the independent 
sequences $((\varepsilon_{1,i}-\varepsilon_{2,i}) \colon i\in\N)$ and 
$((\eta_{1,i,j}-\eta_{2,i,j}) \colon i,j\in\N)$ where each sequence has 
independent and identically distributed elements and satisfies 
Assumption~\ref{as:ts}.
\end{assumption}

\noindent
{\bf Proof of Theorem \ref{thm2}.}
 Theorem~2.1 in \cite{dette2018} implies 
\begin{align*}
\frac{1}{\sqrt{A}}\sum_{i=1}^A (\varepsilon_{1,i} - \varepsilon_{2,i})
\rightsquigarrow Z_\varepsilon \quad \text{and} \quad 
\frac{1}{\sqrt{N}}\sum_{i=1}^A \sum_{j=1}^{n_i} (\eta_{1,i,j}-\eta_{2,i,j})
\rightsquigarrow Z_\eta
\end{align*}
in $C([0,1])$ as $A\to\infty$, $\min_{i=1}^A n_i \to\infty$ where 
$Z_\varepsilon$ and $Z_\eta$ are centred Gaussian processes with covariance kernels 
\begin{align*}
  k_\varepsilon(s,t) &=  \cov (\varepsilon_{1,1}(s)-\varepsilon_{2,1}(s), 
  \varepsilon_{1,1}(t)-\varepsilon_{2,1}(t)) \, , \\
 k_\eta(s,t) &=  \cov (\eta_{1,1,1}(s)-\eta_{2,1,1}(s), 
  \eta_{1,1,1}(t)-\eta_{2,1,1}(t))  \, ,
\end{align*}
respectively. Then we have 
\begin{align} \label{eq:CLT2}
\begin{split}
  \sqrt{A} \, (\hat \theta_N - \theta) 
  &=\frac{1}{\sqrt{A}}\sum_{i=1}^A \big\{ \varepsilon_{1,i} - \varepsilon_{2,i} 
  + \frac{1}{n_i} \sum_{j=1}^{n_i} (\eta_{1,i,j}-\eta_{2,i,j}) \big\} \\
  &= \frac{1}{\sqrt{A}}\sum_{i=1}^A (\varepsilon_{1,i} - \varepsilon_{2,i})
  + o_\P(1) \rightsquigarrow Z_\varepsilon
\end{split}
\end{align} 
in $C([0,1])$ as $A\to\infty$, $\min_{i=1}^A n_i \to\infty$ and similar 
arguments as in the proof of Theorem \ref{thma1} yield
\begin{align} \label{eq:limit-dist3}
\hat{T}^{\theta}_{N} \, {\stackrel{\mathcal{D}}{\longrightarrow}} \,
Z_{\mathcal{E}, \theta} = \max\Big\{ 
\sup_{t\in \mathcal{E}^l_ \kappa} \big(-Z_\varepsilon(t)\big), \,
\sup_{t\in \mathcal{E}^u_ \kappa} Z_\varepsilon(t) \Big\} \, ,
\end{align}
where the statistic $\hat{T}^{\theta}_{N}$ is defined by \eqref{eq:statistic-groups}.

In order to establish the second equality in \eqref{eq:CLT2}, we define 
\begin{align} \label{eq:tilde-eta}
  \tilde \eta_N = \frac{1}{\sqrt{A}}\sum_{i=1}^A 
  \frac{1}{n_i} \sum_{j=1}^{n_i} (\eta_{1,i,j}-\eta_{2,i,j})
\end{align}
and show that this process  converges to zero in 
probability in $C([0,1])$ as $A\to\infty$, $\min_{i=1}^A n_i \to\infty $.
For this purpose we show that the finite dimensional 
distributions converge to zero and prove  that $\tilde \eta_N$ is 
asymptotically $\rho$-equicontinuous in probability, 
where $\rho (s,t) = |s-t|^\theta$. This proves $\tilde \eta_N\rightsquigarrow 0$ in $C([0,1])$  under the stated assumptions
 \citep[see Theorem~7.5 in ][]{billingsley1968}. By the Cram\'er-Wold device, convergence of the finite dimensional distributions to zero is equivalent to 
\begin{align} \label{eq:cramer-wold}
  \sum_{k=1}^q c_k \tilde \eta_N(t_k) 
  = \frac{1}{\sqrt{A}}\sum_{i=1}^A  
  \frac{1}{n_i} \sum_{j=1}^{n_i} \sum_{k=1}^q 
  c_k (\eta_{1,i,j}(t_k) -\eta_{2,i,j}(t_k)) 
  \, {\stackrel{\mathcal{D}}{\longrightarrow}} \, 0
\end{align}
for any $t_1,\dots,t_q\in [0,1]$ and $q\in\N$. Using that the differences $\eta_{1,i,j} - \eta_{2,i,j}$, $i = 1,\dots,A$, $j=1,\dots,n_i$, are independent and the fact that $\E[\eta_{1,i,j}(t_k) -\eta_{2,i,j}(t_k)] = 0$ yields
\begin{align*}
  \E \bigg[ \Big( \sum_{k=1}^q c_j \tilde \eta_N(t_k) \Big)^2 \bigg] 
  &= \frac{1}{A}\sum_{i=1}^A 
  \frac{1}{n_i^2} \sum_{j=1}^{n_i} \E\bigg[ \Big(\sum_{k=1}^q 
  c_k (\eta_{1,i,j}(t_k) -\eta_{2,i,j}(t_k)) \Big)^2\bigg] \\
  &\lesssim \frac{1}{A}\sum_{i=1}^A \frac{1}{n_i}
  \leq \frac{1}{\min_{i=1}^A n_i} \to 0
\end{align*}
where we also used assumption (A1) and the symbol ``$\lesssim$'' means less or equal up to a constant {factor} independent of $A, n_1, \ldots n_A$. This proves \eqref{eq:cramer-wold} and thus the convergence of the finite dimensional distributions to $0$.

In order to verify the equicontinuity condition, we utilize Theorem~2.2.4 in 
\cite{wellner1996}. For any $i=1,\dots,A$, $j=1,\dots,n_i$, define 
$\eta_{i,j} = \eta_{1,i,j} -\eta_{2,i,j}$. Then, for any $s,t\in [0,1]$, we 
have
\begin{align*}
  \E\big[|\tilde \eta_N(s) - \tilde \eta_N(t)|^4\big]^{1/4} \
  &= \E\Big [ \Big | \frac{1}{\sqrt{A}}\sum_{i=1}^A  
  \frac{1}{n_i} \sum_{j=1}^{n_i} (\eta_{i,j}(s) - \eta_{i,j}(t)) 
  \Big |^4 \Big ]^{1/4} \\
  &= \frac{1}{\sqrt{A}} \Big ( \sum_{i=1}^A 
  \frac{1}{n_i^4} \E\Big [ \Big(\sum_{j=1}^{n_i} 
  (\eta_{i,j}(s) - \eta_{i,j}(t)) \Big )^4 \Big ] \\
  &~~ + \sum_{i\neq i'}^A  
  \frac{1}{n_i^2} \frac{1}{n_{i'}^2} \E\Big [ \Big(\sum_{j=1}^{n_i} 
  (\eta_{i,j}(s) - \eta_{i,j}(t)) \Big )^2 \Big ] 
  \E\Big [ \Big( \sum_{j=1}^{n_{i'}} 
  (\eta_{i',j}(s) - \eta_{i',j}(t))  \Big)^2 \Big ] \Big )^{1/4} \, .
\end{align*}
Using assumptions (A1) and (A2) it follows  that 
\begin{align*}
  \E\Big [ \Big(\sum_{j=1}^{n_i}(\eta_{i,j}(s) - \eta_{i,j}(t)) \Big )^4 \Big ] & =
 \sum_{j=1}^{n_i} \E\big[(\eta_{i,j}(s) - \eta_{i,j}(t))^4 \big] \\
  & ~~~
  + \sum_{j\neq j'}^{n_i} \E\big[(\eta_{i,j}(s) - \eta_{i,j}(t))^2 \big]
  \E\big[(\eta_{i,j'}(s) - \eta_{i,j'}(t))^2 \big] \\
  &\lesssim (n_i  + n_i^2) \rho(s,t)^4 \, ,
\\ 
  \E\Big [ \Big(\sum_{j=1}^{n_i} (\eta_{i,j}(s) - \eta_{i,j}(t)) \Big )^2 \Big ] 
 & = \sum_{j=1}^{n_i} \E\big[(\eta_{i,j}(s) - \eta_{i,j}(t))^2 \big]
  \lesssim n_i \rho(s,t)^2 
\end{align*}
which   yields
\begin{align*}
  \E\big[|\tilde \eta_N(s) - \tilde \eta_N(t)|^4\big]^{1/4} 
  \lesssim \frac{1}{\sqrt{A}} \bigg( \sum_{i=1}^A \frac{1}{n_i^2} + \sum_{i\neq i'}^A \frac{1}{n_in_{i'}} \bigg)^{1/4} \rho(s,t) \lesssim \rho(s,t) \, .
\end{align*}
Now we obtain  from Theorem~2.2.4 in \cite{wellner1996} and Markov's inequality that
\begin{align*}
  \P \Big(\sup_{\rho(s,t)\leq \delta} |\tilde \eta_N(s) - \tilde \eta_N(t)| 
  > \varepsilon \Big) 
  &\leq \frac{1}{\varepsilon^4} \E\big[|\tilde\eta_N(s)-\tilde\eta_N(t)|^4\big] \\
  &\lesssim \Big( \int_0^\tau D(\omega,\rho)^{1/4} d\omega + \delta D(\tau,\rho)^{1/2} \Big)^4
\end{align*}
for any $\varepsilon, \delta, \tau > 0$. The discussion {at the beginning of the proof of Theorem~\ref{thma1}} and the fact that $\tau >0$ is arbitrary finally imply
\begin{align*}
  \lim_{\delta \searrow 0} \limsup_{A, \min n_i \to\infty} \P \Big(
  \sup_{\rho(s,t)\leq \delta} |\tilde \eta_N(s) - \tilde \eta_N(t)| 
  > \varepsilon \Big) = 0
\end{align*}
which means that $\tilde\eta_N = o_\P(1)$ since we already proved the convergence of the finite dimensional distributions to zero.

The bootstrap process in \eqref{eq:bootstrap-process2} can be written
\begin{align} \label{eq:boot-epsilon}
\begin{split}
B_N^{\star(r)} &=  \frac{1}{\sqrt{A}}\sum_{i=1}^A 
(M_{i}^{(r)}-1)(\overline{X}_{1,i,\cdot} - \overline{X}_{2,i,\cdot}) \\
&= \frac{1}{\sqrt{A}}\sum_{i=1}^A (M_{i}^{(r)}-1)\big\{ 
\varepsilon_{1,i} - \varepsilon_{2,i} + \frac{1}{n_i} \sum_{j=1}^{n_i} 
(\eta_{1,i,j}-\eta_{2,i,j}) \big\} \\
&= \frac{1}{\sqrt{A}}\sum_{i=1}^A (M_{i}^{(r)}-1)(
\varepsilon_{1,i} - \varepsilon_{2,i}) 
+ o_\P(1)
\end{split}
\end{align}
where $(M_{1}^{(r)},\dots,M_{A}^{(r)})$ follows a multinomial 
distribution with parameters $A$, $p_j = 1/A$, $j = 1,\dots,A$ and the last
estimate follows by  the same arguments as given in the derivation of the second equality in \eqref{eq:CLT2}.

Observe that we have a central limit theorem by the argument given in equation  \eqref{eq:CLT2} and  that   the bootstrap process has a stochastic expansion
given in  \eqref{eq:boot-epsilon}. Therefore it follows by similar arguments as given  in the proof of Theorem~\ref{thm1} that $B_N^{\star(r)}$ converges conditionally given 
$((\varepsilon_{1,i}-\varepsilon_{2,i}) \colon i\in\N)$
 to $Z_\varepsilon$. 
{The continuous mapping theorem and similar arguments as given in the proof of  Lemma~B.3 of \cite{dette2018} yield that the bootstrap statistic defined by \eqref{hd21} converges conditionally given the data $((\varepsilon_{1,i}-\varepsilon_{2,i}) \colon i\in\N)$ to the limit $Z_{\mathcal{E}, \theta}$ in \eqref{eq:limit-dist3}.} 
Furthermore, the assertions in Theorem~\ref{thm2} follow by the same arguments given  in the proof of Theorem~\ref{thm1}.

\subsubsection{Proofs  of the results in Section \ref{sec52} } \label{sec52-assumptions}

\begin{assumption} \label{as:ts2}
  Let $(\eta_{1,i,j}\colon i,j\in\mathbb{N})$ and 
  $(\eta_{2,i,j}\colon i,j\in\mathbb{N})$ denote two (possibly dependent) 
  sequences of $C([0,1])$-valued random variables such that each sequence has 
  independent identically distributed elements and assume that the following 
  conditions are satsified: \vspace{-.2cm}
  \begin{enumerate}\itemsep-.4ex
    \item[(D1)]
    There exist constants $\nu_1, \nu_2>0$, $K$ such that, for all $i,j\in\mathbb{N}$ and $l=1,2$, 
    \begin{align*}
    \mathbb{E}\big[\|\eta_{l,i,j}\|_\infty^{4+\nu_l}\big] \leq K \, , 
    \quad  \mathbb{E}\big[\|\eta_{l,i,j}\|_\infty^{8}\big] < \infty ~.
    \end{align*}
    \item[(D2)]
{There exist constants  $\vartheta>1/4$, $\tilde K$ and  
    real-valued non-negative random variables  $M_1$  and $M_2$ 
    with $\E\big[\|\eta_{l,i,j}\|_\infty^4 M_l^4\big] \leq \tilde K < \infty$} such that, for  $l=1,2$ and  any $i,j\in\mathbb{N}$, 
    the inequality
    \begin{align*}
    |\eta_{l,i,j}(s)-\eta_{l,i,j}(t)|\leq  M_l \, \rho(s,t) = M_l \, |s-t|^\vartheta
    \end{align*}
    holds almost surely for all $s,t\in [0,1]$.
  \end{enumerate}
\end{assumption}

\noindent
{\bf Proof of Theorem \ref{thm3}.}
It is easy to see that Assumption  \ref{as:ts2} implies that the sequences 
of  squared individual errors $(\eta_{1,i,j}^2 \colon i,j\in\N)$ and 
$(\eta_{2,i,j}^2 \colon i,j\in\N)$ satisfy 
Assumption~\ref{as:ts}.
We have for $l=1,2$
\begin{align*}
   \hat{\sigma}_l^2 
   &= \frac{1}{N-A} \sum_{i = 1}^A \sum_{j=1}^{n_i} \Big (X_{l,i,j}
  - \frac{1}{n_i} \sum_{k=1}^{n_i} X_{l,i,k} \Big )^2  \\
  &=  \frac{1}{N-A} \sum_{i = 1}^A \sum_{j=1}^{n_i} \Big (\eta_{l,i,j}
  - \frac{1}{n_i} \sum_{k=1}^{n_i} \eta_{l,i,k} \Big )^2 \\
  &= \frac{1}{N-A} \sum_{i = 1}^A \sum_{j=1}^{n_i} 
   (\eta_{l,i,j})^2 + o_\P \Big ( {1 \over \sqrt{N}}\Big ) ~,
\end{align*}
where the last inequality follows from similar arguments as given in the discussion after equation \eqref{eq:tilde-eta}.
Now we have 
\begin{align*}
  \sqrt{N}&
  \big(\hat{\sigma}_1^2 - \sigma_1^2, \hat{\sigma}_2^2 - \sigma_2^2\big) \\
  &= \sqrt{N} \Big (\frac{1}{N-A} \sum_{i=1}^A \sum_{j=1}^{n_i} 
  \big\{(\eta_{1,i,j})^2 - \sigma_1^2  \big\},
  \frac{1}{N-A} \sum_{i = 1}^A  \sum_{j=1}^{n_i} 
  \big\{(\eta_{2,i,j})^2 - \sigma_2^2  \big\} \Big ) + o_\P(1)
  \rightsquigarrow (Z_1, Z_2)
\end{align*}
in $C([0,1])^2$ as $A\to\infty$, $\min_{i=1}^A n_i \to\infty$. This follows from the fact that  by Theorem~2.1 in 
\cite{dette2018} each component converges individually to its corresponding  limiting process 
in $C([0,1])$. Therefore, both elements are asymptotically 
tight and marginal asymptotic tightness implies joint asymptotic tightness. 
The convergence of the finite dimensional distributions follows from the 
ordinary multidimensional central limit theorem.
By the delta-method as stated in Proposition~2.1 in \cite{carcamo2019} we obtain the convergence 
\begin{align} \label{eq:delta-method}
\begin{split}
  \sqrt{N} &\bigg( \log \bigg(\frac{\hat{\sigma}_1^2}{\hat{\sigma}_2^2} \bigg) - 
  \log \bigg(\frac{\sigma_1^2}{\sigma_2^2}\bigg) \bigg) \\
  &= \frac{\sqrt{N}}{N-A} \sum_{i = 1}^A 
  \sum_{j=1}^{n_i} \Big\{ \frac{(\eta_{1,i,j})^2 - \sigma^2_1}{\sigma^2_1} 
  - \frac{(\eta_{2,i,j})^2 - \sigma^2_2}{\sigma^2_2}  \Big \} + o_\P(1) 
  \rightsquigarrow Z =   \frac{Z_1}{\sigma_1^2} - \frac{Z_2}{\sigma_2^2}
\end{split}
\end{align}
in $ C([0,1]) $, and similar arguments as in the proof of Theorem~\ref{thma1}  yield
\begin{align} \label{hol1}
\hat{T}^{\lambda}_{N} \, {\stackrel{\mathcal{D}}{\longrightarrow}} \,
Z_{\mathcal{E}, \theta} = \max\Big\{ 
\sup_{t\in \mathcal{E}^l_ \zeta} \big(-Z(t)\big), \,
\sup_{t\in \mathcal{E}^u_ \zeta} Z(t) \Big\} \, ,
\end{align}
wh ere the statistic $\hat{T}^{\lambda}_{N}$ is defined by \eqref{eq:statistic-groups_var}.
Turning to the bootstrap process we have for $l=1,2$, $r=1,\dots,R$
\begin{align*}
  C_{l,N}^{\star(r)} &= \frac{1}{N-A} \sum_{i = 1}^A 
  \sum_{j=1}^{n_i} \big\{ (\hat{\eta}_{l,i,j}^{\star (r)})^2 
  - \hat{\sigma}^2_l \big \} \\ 
  &= \frac{1}{N-A} \sum_{i = 1}^A 
  \sum_{j=1}^{n_i} \big\{ (\hat{\eta}_{l,i,j}^{\star (r)})^2 
  - (\hat{\eta}_{l,i,j})^2 \big \} + {o_\P(1)} \\
  &=\frac{1}{N-A} \sum_{i = 1}^A 
  \sum_{j=1}^{n_i} (M_{i,j}^{(r)}-1) 
  \Big (\eta_{l,i,j}
  - \frac{1}{n_i} \sum_{k=1}^{n_i} \eta_{l,i,k} \Big )^2 + {o_\P(1)} \\
  &=\frac{1}{N-A} \sum_{i = 1}^A 
  \sum_{j=1}^{n_i} (M_{i,j}^{(r)}-1) (\eta_{l,i,j})^2 
  + o_\P \Big ({1 \over \sqrt{N}}  \Big ) \, ,
\end{align*}
where the last estimate follows from similar arguments as those used to derive  the second equality in \eqref{eq:CLT2} 
and the vector $(M_{1,1}^{(r)},\ldots,M_{1,n_1}^{(r)},\ldots , M_{A,1}^{(r)}, \ldots,  M_{A,n_A}^{(r)})$ follows a multinomial 
  distribution with parameters $N$ and equal probabilities $ 1/N$.
Then 
\begin{align*}
 \sqrt{N} \, \tilde C_{N}^{\star(r)} &= \sqrt{N} \Big ( 
 \frac{C_{1,N}^{\star(r)}}{\sigma^2_1} - \frac{C_{2,N}^{\star(r)}}{\sigma^2_2} 
 \Big ) \\
 &= \frac{\sqrt{N}}{N-A} \sum_{i = 1}^A \sum_{j=1}^{n_i} (M_{i,j}^{(r)}-1) 
 \Big \{\frac{(\eta_{1,i,j})^2}{\sigma_1^2} - \frac{(\eta_{2,i,j})^2}{\sigma_2^2} \Big \} + o_\P(1)
\end{align*}
and this process converges conditionally given {$(\eta_{1,i,j}\colon i,j\in\mathbb{N})$ and 
  $(\eta_{2,i,j}\colon i,j\in\mathbb{N})$}  to the same limit as the limit in \eqref{eq:delta-method}. Finally it can be shown that
\begin{align*}
  \sqrt{N}\|C_{N}^{\star(r)} - \tilde C_{N}^{\star(r)} \|_\infty 
  = \sqrt{N}\|C_{1,N}^{\star(r)}(1/\sigma_1^2 - 1/\hat\sigma_1^2) 
  - C_{2,N}^{\star(r)} (1/\sigma_2^2 - 1/\hat\sigma_2^2)\|_\infty
  = o_\P(1)
\end{align*}
and by the continuous mapping theorem and similar arguments as given in the proof of  Lemma~B.3 of  \cite{dette2018}, the bootstrap statistic defined by \eqref{eq:bootstrap-statistic-var} converges conditionally given the data $(X_{lj} \colon j\in \N; ~ l=1,2)$ to the  limit in \eqref{hol1}.

The assertion of Theorem~\ref{thm3} now follows by the same arguments given in the proof of {Theorem~\ref{thm1}}.

\end{document}